\documentclass[titlepage,a4paper,12pt]{article}
\usepackage{latexsym,amssymb}
\font\Fraktur=eufm10 scaled\magstep0          
\newcommand{\fraktur}[1]{\mbox{\Fraktur #1}}  %
\font\Fraktu=eufm7 scaled\magstep0            
\newcommand{\fraktu}[1]{\mbox{\Fraktu #1}}    %
\font\Frakt=eufm5 scaled\magstep0             
\newcommand{\frakt}[1]{\mbox{\Frakt #1}}      %
\def\fr#1{\mathchoice{\fraktur {#1}}          
                     {\fraktur {#1}}          
                     {\fraktu {#1}}           
                     {\frakt {#1}}}           
\newtheorem{Definition}{Definition}[section]
\newtheorem{Theorem and Definition}[Definition]{Theorem and Definition}
\newtheorem{Proposition and Definition}[Definition]{Proposition and Definition}
\newtheorem{Example}[Definition]{Example}
\newtheorem{Theorem}[Definition]{Theorem}
\newtheorem{Proposition}[Definition]{Proposition}
\newtheorem{Remark}[Definition]{Remark}

\newtheorem{Corollary}[Definition]{Corollary}

\newtheorem{Note}[Definition]{Note}
\newtheorem{Par}[Definition]{} 
\newenvironment{Proof}{{\noindent \sc Proof: } }{\mbox{ }\hfill$\Box$  
  \vspace{1.5ex} \par} 
\newenvironment{rPar}{\begin{Par} \em}{ \end{Par}}
\newenvironment{rExample}{\begin{Example} \em}{ \end{Example}}
\newenvironment{rRemark}{\begin{Remark} \em}{ \end{Remark}}

\newcommand{\N}{{\Bbb{N}}}
\newcommand{\C}{{\Bbb{C}}}
\newcommand{\R}{{\Bbb{R}}}
\newcommand{\K}{{\Bbb{K}}}
\newcommand{\Q}{{\Bbb{Q}}}
\newcommand{\Z}{{\Bbb{Z}}}

\newcommand{\m}{\fr m}

\newcommand{\g}{\fr g}

\newcommand{\e}{{\rm e}}
\newcommand{\id}{{{\rm id}}}
\newcommand{\GL}{{\rm GL \,}}
\newcommand{\End}{{\rm End \,}}

\newcommand{\cat}[1]{\underline{\mbox{\sf {#1}}}}

\newcommand{\dera}[1]{\frac{\partial^{|\alpha|}}{\partial {#1}^{\alpha}}}
\newcommand{\derb}[1]{\frac{\partial^{|\beta|}}{\partial {#1}^{\beta}}}

\newcommand{\im}[1]{{\rm im} \: {#1} }

\newcommand{\symi}{{\rm S}_{\rho,\delta}^{\infty}}

\newcommand{\symmi}{{\rm S}^{-\infty}}

\newcommand{\hattimes}{\hat{\otimes}}
\long\def\ig#1{\relax}
\ig{Thanks to Roberto Minio for this def'n.  Compare the def'n of
\comment in AMSTeX.}
 
\newcount \coefa
\newcount \coefb
\newcount \coefc
\newdimen\tempdimen
\newdimen\xlen
\newdimen\ylen
 
\newcount\tempcounta
\newcount\tempcountb
\newcount\tempcountc
\newcount\tempcountd
\newcount\tempcounte
\newcount\tempcountf
\newcount\xext
\newcount\yext
\newcount\xoff
\newcount\yoff
\newcount\gap%
\newcount\arrowtypea
\newcount\arrowtypeb
\newcount\arrowtypec
\newcount\arrowtyped
\newcount\arrowtypee
\newcount\height
\newcount\width
\newcount\xpos
\newcount\ypos
\newcount\run
\newcount\rise
\newcount\arrowlength
\newcount\halflength
\newcount\arrowtype
\newsavebox{\tempboxa}%
\newsavebox{\tempboxb}%
\newsavebox{\tempboxc}%

\catcode`@=11 
\def\settoheight#1#2{\setbox\@tempboxa\hbox{#2}#1\ht\@tempboxa\relax}%
\def\settodepth#1#2{\setbox\@tempboxa\hbox{#2}#1\dp\@tempboxa\relax}%
\let\ifnextchar=\@ifnextchar
\catcode`@=12 
 
\def\putbox(#1,#2)#3{\put(#1,#2){\makebox(0,0){#3}}}

\def\setsqparms[#1`#2`#3`#4;#5`#6]{%
\settripairparms[#1`#2`#3`#4`1;#6]%
\width #5
}
 
\def\settriparms[#1`#2`#3;#4]{\settripairparms[#1`#2`#3`1`1;#4]}%

\def\settripairparms[#1`#2`#3`#4`#5;#6]{%
\arrowtypea #1
\arrowtypeb #2
\arrowtypec #3
\arrowtyped #4
\arrowtypee #5
\height #6
\width #6
}
 
\def\mvector(#1,#2)#3{
\put(0,0){\vector(#1,#2){#3}}%
\put(0,0){\vector(#1,#2){30}}%
}
\def\evector(#1,#2)#3{{
\arrowlength #3
\put(0,0){\vector(#1,#2){\arrowlength}}%
\advance \arrowlength by-30
\put(0,0){\vector(#1,#2){\arrowlength}}%
}}

\def\horsize#1#2{%
\settowidth{\tempdimen}{$#2$}%
#1=\tempdimen
\divide #1 by\unitlength
}
 
\def\vertsize#1#2{%
\settoheight{\tempdimen}{$#2$}%
#1=\tempdimen
\settodepth{\tempdimen}{$#2$}%
\advance #1 by\tempdimen
\divide #1 by\unitlength
}

\def\vertadjust[#1`#2`#3]{%
\vertsize{\tempcounta}{#1}%
\vertsize{\tempcountb}{#2}%
\ifnum \tempcounta<\tempcountb \tempcounta=\tempcountb \fi
\divide\tempcounta by2
\vertsize{\tempcountb}{#3}%
\ifnum \tempcountb>0 \advance \tempcountb by20 \fi
\ifnum \tempcounta<\tempcountb \tempcounta=\tempcountb \fi
}
 
\def\horadjust[#1`#2`#3]{%
\horsize{\tempcounta}{#1}%
\horsize{\tempcountb}{#2}%
\ifnum \tempcounta<\tempcountb \tempcounta=\tempcountb \fi
\divide\tempcounta by20
\horsize{\tempcountb}{#3}%
\ifnum \tempcountb>0 \advance \tempcountb by60 \fi
\ifnum \tempcounta<\tempcountb \tempcounta=\tempcountb \fi
}
 
\ig{ In this procedure, #1 is the paramater that sticks out all the way,
#2 sticks out the least and #3 is a label sticking out half way.  #4 is
the amount of the offset.}
 
\def\sladjust[#1`#2`#3]#4{%
\tempcountc=#4
\horsize{\tempcounta}{#1}%
\divide \tempcounta by2
\horsize{\tempcountb}{#2}%
\divide \tempcountb by2
\advance \tempcountb by-\tempcountc
\ifnum \tempcounta<\tempcountb \tempcounta=\tempcountb\fi
\divide \tempcountc by2
\horsize{\tempcountb}{#3}%
\advance \tempcountb by-\tempcountc
\ifnum \tempcountb>0 \advance \tempcountb by80\fi
\ifnum \tempcounta<\tempcountb \tempcounta=\tempcountb\fi
\advance\tempcounta by20
}
 
\def\putvector(#1,#2)(#3,#4)#5#6{{%
\xpos=#1
\ypos=#2
\run=#3
\rise=#4
\arrowlength=#5
\arrowtype=#6
\ifnum \arrowtype<0
    \ifnum \run=0
        \advance \ypos by-\arrowlength
    \else
        \tempcounta \arrowlength
        \multiply \tempcounta by\rise
        \divide \tempcounta by\run
        \ifnum\run>0
            \advance \xpos by\arrowlength
            \advance \ypos by\tempcounta
        \else
            \advance \xpos by-\arrowlength
            \advance \ypos by-\tempcounta
        \fi
    \fi
    \multiply \arrowtype by-1
    \multiply \rise by-1
    \multiply \run by-1
\fi
\ifnum \arrowtype=1
    \put(\xpos,\ypos){\vector(\run,\rise){\arrowlength}}%
\else\ifnum \arrowtype=2
    \put(\xpos,\ypos){\mvector(\run,\rise)\arrowlength}%
\else\ifnum\arrowtype=3
    \put(\xpos,\ypos){\evector(\run,\rise){\arrowlength}}%
\fi\fi\fi
}}
 
\def\bfig{\begin{picture}(\xext,\yext)(\xoff,\yoff)}
\def\efig{\end{picture}}

\def\putsplitvector(#1,#2)#3#4{
\xpos #1
\ypos #2
\arrowtype #4
\halflength #3
\arrowlength #3
\gap 140
\advance \halflength by-\gap
\divide \halflength by2
\ifnum \arrowtype=1
    \put(\xpos,\ypos){\line(0,-1){\halflength}}%
    \advance\ypos by-\halflength
    \advance\ypos by-\gap
    \put(\xpos,\ypos){\vector(0,-1){\halflength}}%
\else\ifnum \arrowtype=2
    \put(\xpos,\ypos){\line(0,-1)\halflength}%
    \put(\xpos,\ypos){\vector(0,-1)3}%
    \advance\ypos by-\halflength
    \advance\ypos by-\gap
    \put(\xpos,\ypos){\vector(0,-1){\halflength}}%
\else\ifnum\arrowtype=3
    \put(\xpos,\ypos){\line(0,-1)\halflength}%
    \advance\ypos by-\halflength
    \advance\ypos by-\gap
    \put(\xpos,\ypos){\evector(0,-1){\halflength}}%
\else\ifnum \arrowtype=-1
    \advance \ypos by-\arrowlength
    \put(\xpos,\ypos){\line(0,1){\halflength}}%
    \advance\ypos by\halflength
    \advance\ypos by\gap
    \put(\xpos,\ypos){\vector(0,1){\halflength}}%
\else\ifnum \arrowtype=-2
    \advance \ypos by-\arrowlength
    \put(\xpos,\ypos){\line(0,1)\halflength}%
    \put(\xpos,\ypos){\vector(0,1)3}%
    \advance\ypos by\halflength
    \advance\ypos by\gap
    \put(\xpos,\ypos){\vector(0,1){\halflength}}%
\else\ifnum\arrowtype=-3
    \advance \ypos by-\arrowlength
    \put(\xpos,\ypos){\line(0,1)\halflength}%
    \advance\ypos by\halflength
    \advance\ypos by\gap
    \put(\xpos,\ypos){\evector(0,1){\halflength}}%
\fi\fi\fi\fi\fi\fi
}
 
\def\setpos(#1,#2){\xpos=#1 \ypos#2}
 
\def\putmorphism(#1)(#2,#3)[#4`#5`#6]#7#8#9{{%
\run #2
\rise #3
\ifnum\rise=0
  \puthmorphism(#1)[#4`#5`#6]{#7}{#8}{#9}%
\else\ifnum\run=0
  \putvmorphism(#1)[#4`#5`#6]{#7}{#8}{#9}%
\else
\setpos(#1)%
\arrowlength #7
\arrowtype #8
\ifnum\run=0
\else\ifnum\rise=0
\else
\ifnum\run>0
    \coefa=1
\else
   \coefa=-1
\fi
\ifnum\arrowtype>0
   \coefb=0
   \coefc=-1
\else
   \coefb=\coefa
   \coefc=1
   \arrowtype=-\arrowtype
\fi
\width=2
\multiply \width by\run
\divide \width by\rise
\ifnum \width<0  \width=-\width\fi
\advance\width by60
\if l#9 \width=-\width\fi
\putbox(\xpos,\ypos){$#4$}
{\multiply \coefa by\arrowlength
\advance\xpos by\coefa
\multiply \coefa by\rise
\divide \coefa by\run
\advance \ypos by\coefa
\putbox(\xpos,\ypos){$#5$} }%
{\multiply \coefa by\arrowlength
\divide \coefa by2
\advance \xpos by\coefa
\advance \xpos by\width
\multiply \coefa by\rise
\divide \coefa by\run
\advance \ypos by\coefa
\if l#9%
   \put(\xpos,\ypos){\makebox(0,0)[r]{$#6$}}%
\else\if r#9%
   \put(\xpos,\ypos){\makebox(0,0)[l]{$#6$}}%
\fi\fi }%
{\multiply \rise by-\coefc
\multiply \run by-\coefc
\multiply \coefb by\arrowlength
\advance \xpos by\coefb
\multiply \coefb by\rise
\divide \coefb by\run
\advance \ypos by\coefb
\multiply \coefc by70
\advance \ypos by\coefc
\multiply \coefc by\run
\divide \coefc by\rise
\advance \xpos by\coefc
\multiply \coefa by140
\multiply \coefa by\run
\divide \coefa by\rise
\advance \arrowlength by\coefa
\ifnum \arrowtype=1
   \put(\xpos,\ypos){\vector(\run,\rise){\arrowlength}}%
\else\ifnum\arrowtype=2
   \put(\xpos,\ypos){\mvector(\run,\rise){\arrowlength}}%
\else\ifnum\arrowtype=3
   \put(\xpos,\ypos){\evector(\run,\rise){\arrowlength}}%
\fi\fi\fi}%
\fi\fi
\fi\fi}}

\def\puthmorphism(#1,#2)[#3`#4`#5]#6#7#8{{%
\xpos #1
\ypos #2
\width #6
\arrowlength #6
\putbox(\xpos,\ypos){$#3$\vphantom{$#4$}}%
{\advance \xpos by\arrowlength
\putbox(\xpos,\ypos){\vphantom{$#3$}$#4$}}%
\horsize{\tempcounta}{#3}%
\horsize{\tempcountb}{#4}%
\divide \tempcounta by2
\divide \tempcountb by2
\advance \tempcounta by30
\advance \tempcountb by30
\advance \xpos by\tempcounta
\advance \arrowlength by-\tempcounta
\advance \arrowlength by-\tempcountb
\putvector(\xpos,\ypos)(1,0){\arrowlength}{#7}%
\divide \arrowlength by2
\advance \xpos by\arrowlength
\vertsize{\tempcounta}{#5}%
\divide\tempcounta by2
\advance \tempcounta by20
\if a#8 %
   \advance \ypos by\tempcounta
   \put(\xpos,\ypos){\makebox(0,0){$#5$}}%
\else
   \advance \ypos by-\tempcounta
   \put(\xpos,\ypos){\makebox(0,0){$#5$}}%
\fi
}}
 
\def\putvmorphism(#1,#2)[#3`#4`#5]#6#7#8{{%
\xpos #1
\ypos #2
\arrowlength #6
\arrowtype #7
\settowidth{\xlen}{$#5$}%
\putbox(\xpos,\ypos){$#3$}%
{\advance \ypos by-\arrowlength
\putbox(\xpos,\ypos){$#4$}}%
{\advance\arrowlength by-140
\advance \ypos by-70
\ifdim\xlen>0pt
   \if m#8%
      \putsplitvector(\xpos,\ypos){\arrowlength}{\arrowtype}%
   \else
      \putvector(\xpos,\ypos)(0,-1){\arrowlength}{\arrowtype}%
   \fi
\else
   \putvector(\xpos,\ypos)(0,-1){\arrowlength}{\arrowtype}%
\fi}%
\ifdim\xlen>0pt
   \divide \arrowlength by2
   \advance\ypos by-\arrowlength
   \if l#8%
      \advance \xpos by-40
      \put(\xpos,\ypos){\makebox(0,0)[r]{$#5$}}%
   \else\if r#8%
      \advance \xpos by40
      \put(\xpos,\ypos){\makebox(0,0)[l]{$#5$}}%
   \else
      \putbox(\xpos,\ypos){$#5$}%
   \fi\fi
\fi
}}
 
\def\topadjust[#1`#2`#3]{%
\yoff=10
\vertadjust[#1`#2`{#3}]%
\advance \yext by\tempcounta
\advance \yext by 10
}
 
\def\botadjust[#1`#2`#3]{%
\vertadjust[#1`#2`{#3}]%
\advance \yext by\tempcounta
\advance \yoff by-\tempcounta
}
 
\def\leftadjust[#1`#2`#3]{%
\xoff=0
\horadjust[#1`#2`{#3}]%
\advance \xext by\tempcounta
\advance \xoff by-\tempcounta
}
 
\def\rightadjust[#1`#2`#3]{%
\horadjust[#1`#2`{#3}]%
\advance \xext by\tempcounta
}
 
\def\rightsladjust[#1`#2`#3]{%
\sladjust[#1`#2`{#3}]{\width}%
\advance \xext by\tempcounta
}
 
\def\leftsladjust[#1`#2`#3]{%
\xoff=0
\sladjust[#1`#2`{#3}]{\width}%
\advance \xext by\tempcounta
\advance \xoff by-\tempcounta
}
 
\def\adjust[#1`#2;#3`#4;#5`#6;#7`#8]{%
\topadjust[#1``{#2}]
\leftadjust[#3``{#4}]
\rightadjust[#5``{#6}]
\botadjust[#7``{#8}]}

\def\putsquare(#1)[#2`#3`#4`#5;#6`#7`#8`#9]{%
\setpos(#1)
\puthmorphism(\xpos,\ypos)[#4`#5`{#9}]{\width}{\arrowtyped}b%
\advance\ypos by \height
\puthmorphism(\xpos,\ypos)[#2`#3`{#6}]{\width}{\arrowtypea}a%
\putvmorphism(\xpos,\ypos)[``{#7}]{\height}{\arrowtypeb}l%
\advance\xpos by \width
\putvmorphism(\xpos,\ypos)[``{#8}]{\height}{\arrowtypec}r%
}
 
\def\square[#1`#2`#3`#4;#5`#6`#7`#8]{{
\xext=\width                              
\yext=\height                             
\topadjust[#1`#2`{#5}]
\botadjust[#3`#4`{#8}]
\leftadjust[#1`#3`{#6}]
\rightadjust[#2`#4`{#7}]
\begin{picture}(\xext,\yext)(\xoff,\yoff)
\putsquare(0,0)[#1`#2`#3`#4;#5`#6`#7`{#8}]
\end{picture}
}}

\def\putptriangle(#1,#2)[#3`#4`#5;#6`#7`#8]{%
\xpos=#1 \ypos=#2
\advance\ypos by \height
\puthmorphism(\xpos,\ypos)[#3`#4`{#6}]{\height}{\arrowtypea}a%
\putvmorphism(\xpos,\ypos)[`#5`{#7}]{\height}{\arrowtypeb}l%
\advance\xpos by\height
\putmorphism(\xpos,\ypos)(-1,-1)[``{#8}]{\height}{\arrowtypec}r%
}
 
\def\ptriangle[#1`#2`#3;#4`#5`#6]{{
\width=\height                         
\xext=\width                           
\yext=\width                           
\topadjust[#1`#2`{#4}]
\botadjust[#3``]
\leftadjust[#1`#3`{#5}]
\rightsladjust[#2`#3`{#6}]
\begin{picture}(\xext,\yext)(\xoff,\yoff)
\putptriangle(0,0)[#1`#2`#3;#4`#5`{#6}]%
\end{picture}%
}}

\def\putqtriangle(#1,#2)[#3`#4`#5;#6`#7`#8]{%
\xpos=#1 \ypos=#2
\advance\ypos by\height
\puthmorphism(\xpos,\ypos)[#3`#4`{#6}]{\height}{\arrowtypea}a%
\putmorphism(\xpos,\ypos)(1,-1)[``{#7}]{\height}{\arrowtypeb}l%
\advance\xpos by\height
\putvmorphism(\xpos,\ypos)[`#5`{#8}]{\height}{\arrowtypec}r%
}
 
\def\qtriangle[#1`#2`#3;#4`#5`#6]{{
\width=\height                         
\xext=\width                           
\yext=\height                          
\topadjust[#1`#2`{#4}]
\botadjust[#3``]
\leftsladjust[#1`#3`{#5}]
\rightadjust[#2`#3`{#6}]
\begin{picture}(\xext,\yext)(\xoff,\yoff)
\putqtriangle(0,0)[#1`#2`#3;#4`#5`{#6}]%
\end{picture}%
}}

\def\putdtriangle(#1,#2)[#3`#4`#5;#6`#7`#8]{%
\xpos=#1 \ypos=#2
\puthmorphism(\xpos,\ypos)[#4`#5`{#8}]{\height}{\arrowtypec}b%
\advance\xpos by \height \advance\ypos by\height
\putmorphism(\xpos,\ypos)(-1,-1)[``{#6}]{\height}{\arrowtypea}l%
\putvmorphism(\xpos,\ypos)[#3``{#7}]{\height}{\arrowtypeb}r%
}
 
\def\dtriangle[#1`#2`#3;#4`#5`#6]{{
\width=\height                         
\xext=\width                           
\yext=\height                          
\topadjust[#1``]
\botadjust[#2`#3`{#6}]
\leftsladjust[#2`#1`{#4}]
\rightadjust[#1`#3`{#5}]
\begin{picture}(\xext,\yext)(\xoff,\yoff)
\putdtriangle(0,0)[#1`#2`#3;#4`#5`{#6}]%
\end{picture}%
}}

\def\putbtriangle(#1,#2)[#3`#4`#5;#6`#7`#8]{%
\xpos=#1 \ypos=#2
\puthmorphism(\xpos,\ypos)[#4`#5`{#8}]{\height}{\arrowtypec}b%
\advance\ypos by\height
\putmorphism(\xpos,\ypos)(1,-1)[``{#7}]{\height}{\arrowtypeb}r%
\putvmorphism(\xpos,\ypos)[#3``{#6}]{\height}{\arrowtypea}l%
}
 
\def\btriangle[#1`#2`#3;#4`#5`#6]{{
\width=\height                         
\xext=\width                           
\yext=\height                          
\topadjust[#1``]
\botadjust[#2`#3`{#6}]
\leftadjust[#1`#2`{#4}]
\rightsladjust[#3`#1`{#5}]
\begin{picture}(\xext,\yext)(\xoff,\yoff)
\putbtriangle(0,0)[#1`#2`#3;#4`#5`{#6}]%
\end{picture}%
}}

\def\putAtriangle(#1,#2)[#3`#4`#5;#6`#7`#8]{%
\xpos=#1 \ypos=#2
{\multiply \height by2
\puthmorphism(\xpos,\ypos)[#4`#5`{#8}]{\height}{\arrowtypec}b}%
\advance\xpos by\height \advance\ypos by\height
\putmorphism(\xpos,\ypos)(-1,-1)[#3``{#6}]{\height}{\arrowtypea}l%
\putmorphism(\xpos,\ypos)(1,-1)[``{#7}]{\height}{\arrowtypeb}r%
}
 
\def\Atriangle[#1`#2`#3;#4`#5`#6]{{
\width=\height                         
\xext=\width                           
\yext=\height                          
\topadjust[#1``]
\botadjust[#2`#3`{#6}]
\multiply \xext by2 
\leftsladjust[#2`#1`{#4}]
\rightsladjust[#3`#1`{#5}]
\begin{picture}(\xext,\yext)(\xoff,\yoff)%
\putAtriangle(0,0)[#1`#2`#3;#4`#5`{#6}]%
\end{picture}%
}}

\def\putAtrianglepair(#1,#2)[#3]{\xpos=#1 \ypos=#2%
\putAtrianglepairx[#3]}
\def\putAtrianglepairx[#1`#2`#3`#4;#5`#6`#7`#8`#9]{%
\puthmorphism(\xpos,\ypos)[#2`#3`{#8}]{\height}{\arrowtyped}b%
\advance\xpos by\height
\puthmorphism(\xpos,\ypos)[\phantom{#3}`#4`{#9}]{\height}{\arrowtypee}b%
\advance\ypos by\height
\putmorphism(\xpos,\ypos)(-1,-1)[#1``{#5}]{\height}{\arrowtypea}l%
\putvmorphism(\xpos,\ypos)[``{#6}]{\height}{\arrowtypeb}m%
\putmorphism(\xpos,\ypos)(1,-1)[``{#7}]{\height}{\arrowtypec}r%
}
 
\def\Atrianglepair[#1`#2`#3`#4;#5`#6`#7`#8`#9]{{%
\width=\height
\xext=\width
\yext=\height
\topadjust[#1``]%
\vertadjust[#2`#3`{#8}]
\tempcountd=\tempcounta                       
\vertadjust[#3`#4`{#9}]
\ifnum\tempcounta<\tempcountd                 
\tempcounta=\tempcountd\fi                    
\advance \yext by\tempcounta                  
\advance \yoff by-\tempcounta                 
\multiply \xext by2 
\leftsladjust[#2`#1`{#5}]
\rightsladjust[#4`#1`{#7}]%
\begin{picture}(\xext,\yext)(\xoff,\yoff)%
\putAtrianglepair(0,0)[#1`#2`#3`#4;#5`#6`#7`#8`{#9}]%
\end{picture}%
}}

\def\putVtriangle(#1,#2)[#3`#4`#5;#6`#7`#8]{%
\xpos=#1 \ypos=#2
\advance\ypos by\height
{\multiply\height by2
\puthmorphism(\xpos,\ypos)[#3`#4`{#6}]{\height}{\arrowtypea}a}%
\putmorphism(\xpos,\ypos)(1,-1)[`#5`{#7}]{\height}{\arrowtypeb}l%
\advance\xpos by\height
\advance\xpos by\height
\putmorphism(\xpos,\ypos)(-1,-1)[``{#8}]{\height}{\arrowtypec}r%
}
 
\def\Vtriangle[#1`#2`#3;#4`#5`#6]{{
\width=\height                         
\xext=\width                           
\yext=\height                          
\topadjust[#1`#2`{#4}]
\botadjust[#3``]
\multiply \xext by2 
\leftsladjust[#1`#3`{#5}]
\rightsladjust[#2`#3`{#6}]
\begin{picture}(\xext,\yext)(\xoff,\yoff)%
\putVtriangle(0,0)[#1`#2`#3;#4`#5`{#6}]%
\end{picture}%
}}

\def\putVtrianglepair(#1,#2)[#3]{\xpos=#1 \ypos=#2%
\putVtrianglepairx[#3]}
\def\putVtrianglepairx[#1`#2`#3`#4;#5`#6`#7`#8`#9]{%
\advance\ypos by\height
\putmorphism(\xpos,\ypos)(1,-1)[`#4`{#7}]{\height}{\arrowtypec}l%
\puthmorphism(\xpos,\ypos)[#1`#2`{#5}]{\height}{\arrowtypea}a%
\advance\xpos by\height
\puthmorphism(\xpos,\ypos)[\phantom{#2}`#3`{#6}]{\height}{\arrowtypeb}a%
\putvmorphism(\xpos,\ypos)[``{#8}]{\height}{\arrowtyped}m%
\advance\xpos by\height
\putmorphism(\xpos,\ypos)(-1,-1)[``{#9}]{\height}{\arrowtypee}r%
}
 
\def\Vtrianglepair[#1`#2`#3`#4;#5`#6`#7`#8`#9]{{%
\xoff=0
\yoff=2 
\xext=\height                  
\width=\height                 
\yext=\height                  
\vertadjust[#1`#2`{#5}]
\tempcountd=\tempcounta        
\vertadjust[#2`#3`{#6}]
\ifnum\tempcounta<\tempcountd
\tempcounta=\tempcountd\fi
\advance \yext by\tempcounta
\botadjust[#4``]%
\multiply \xext by2
\leftsladjust[#1`#4`{#7}]%
\rightsladjust[#3`#4`{#9}]%
\begin{picture}(\xext,\yext)(\xoff,\yoff)%
\putVtrianglepair(0,0)[#1`#2`#3`#4;#5`#6`#7`#8`{#9}]%
\end{picture}%
}}

\def\putCtriangle(#1,#2)[#3`#4`#5;#6`#7`#8]{%
\xpos=#1 \ypos=#2
\advance\ypos by\height
\putmorphism(\xpos,\ypos)(1,-1)[``{#8}]{\height}{\arrowtypec}l%
\advance\xpos by\height
\advance\ypos by\height
\putmorphism(\xpos,\ypos)(-1,-1)[#3`#4`{#6}]{\height}{\arrowtypea}l%
{\multiply\height by 2
\putvmorphism(\xpos,\ypos)[`#5`{#7}]{\height}{\arrowtypeb}r}%
}
 
\def\Ctriangle[#1`#2`#3;#4`#5`#6]{{
\width=\height                          
\xext=\width                            
\yext=\height                           
\multiply \yext by2 
\topadjust[#1``]
\botadjust[#3``]
\sladjust[#2`#1`{#4}]{\width}
\tempcountd=\tempcounta                 
\sladjust[#2`#3`{#6}]{\width}
\ifnum \tempcounta<\tempcountd          
\tempcounta=\tempcountd\fi              
\advance \xext by\tempcounta            
\advance \xoff by-\tempcounta           
\rightadjust[#1`#3`{#5}]
\begin{picture}(\xext,\yext)(\xoff,\yoff)%
\putCtriangle(0,0)[#1`#2`#3;#4`#5`{#6}]%
\end{picture}%
}}

\def\putDtriangle(#1,#2)[#3`#4`#5;#6`#7`#8]{%
\xpos=#1 \ypos=#2
\advance\xpos by\height \advance\ypos by\height
\putmorphism(\xpos,\ypos)(-1,-1)[``{#8}]{\height}{\arrowtypec}r%
\advance\xpos by-\height \advance\ypos by\height
\putmorphism(\xpos,\ypos)(1,-1)[`#4`{#7}]{\height}{\arrowtypeb}r%
{\multiply\height by 2
\putvmorphism(\xpos,\ypos)[#3`#5`{#6}]{\height}{\arrowtypea}l}%
}
 
\def\Dtriangle[#1`#2`#3;#4`#5`#6]{{
\width=\height                         
\xext=\height                          
\yext=\height                          
\multiply \yext by2 
\topadjust[#1``]
\botadjust[#3``]
\leftadjust[#1`#3`{#4}]
\sladjust[#2`#1`{#4}]{\height}
\tempcountd=\tempcountd                
\sladjust[#2`#3`{#6}]{\height}
\ifnum \tempcounta<\tempcountd         
\tempcounta=\tempcountd\fi             
\advance \xext by\tempcounta           
\begin{picture}(\xext,\yext)(\xoff,\yoff)
\putDtriangle(0,0)[#1`#2`#3;#4`#5`{#6}]%
\end{picture}%
}}

\def\setrecparms[#1`#2]{\width=#1 \height=#2}%
\def\recurse[#1`#2`#3`#4;#5`#6`#7`#8`#9]{{%
\settowidth{\tempdimen}{#1}
\ifdim\tempdimen=0pt
  \savebox{\tempboxa}{\hbox{#2}}%
  \savebox{\tempboxb}{\hbox{#4}}%
  \savebox{\tempboxc}{\hbox{#7}}%
\else
  \savebox{\tempboxa}{\hbox{$\hbox{#1}\times\hbox{#2}$}}%
  \savebox{\tempboxb}{\hbox{$\hbox{#1}\times\hbox{#4}$}}%
  \savebox{\tempboxc}{\hbox{$\hbox{#1}\times\hbox{#7}$}}%
\fi
\tempcounte=\height
\divide\tempcounte by 2
\tempcountf=\tempcounte
\advance\tempcountf by \width
\xext=\tempcountf \yext=\height
\topadjust[#2`\usebox{\tempboxa}`{#5}]%
\botadjust[#4`\usebox{\tempboxb}`{#9}]%
\sladjust[#3`#2`{#6}]{\tempcounte}%
\tempcountd=\tempcounta
\sladjust[#3`#4`{#8}]{\tempcounte}%
\ifnum \tempcounta<\tempcountd
\tempcounta=\tempcountd\fi
\advance \xext by\tempcounta
\advance \xoff by-\tempcounta
\rightadjust[\usebox{\tempboxa}`\usebox{\tempboxb}`\usebox{\tempboxc}]%
\bfig
{\settriparms[-1`1`1;\tempcounte]%
\putCtriangle(0,0)[`#3`;#6`#7`{#8}]}%
\arrowtypea=-1 \arrowtypeb=0 \arrowtypec=1 \arrowtyped=-1
\putsquare(\tempcounte,0)[#2`\usebox{\tempboxa}`#4`\usebox{\tempboxb};%
#5``\usebox{\tempboxc}`#9]%
\efig
}}

\begin{document}
\begin{titlepage}
\thispagestyle{empty}
\title{Holomorphic deformation of Hopf algebras and applications to quantum
       groups}
\author{ Markus J.~Pflaum\thanks{ Humboldt Universit\"at zu Berlin, 
         Mathematisches Institut, Unter den Linden 6, 10099 Berlin,
         \mbox{ } \hspace{1cm} email: pflaum@mathematik.hu-berlin.de }\\
         Martin Schottenloher\thanks{\mbox{Ludwig-Maximilians-Universit\"at, 
         Mathematisches Institut, Theresienstr.~39, 80333 M\"unchen,}
	\protect\newline
         \mbox{ } \hspace{1cm} email: schotten@rz.mathematik.uni-muenchen.de}}    
\date{ \today }
\maketitle
\end{titlepage}
\begin{abstract}
In this article we propose a new and so-called holomorphic deformation scheme 
for locally convex algebras and Hopf algebras. 
Essentially we regard converging power series expansion of a deformed product 
on a locally convex algebra, thus giving the means to actually insert 
complex values for the deformation parameter.  
Moreover we establish a topological duality theory for locally convex
Hopf algebras.
Examples coming from the theory of quantum groups are reconsidered 
within our holomorphic deformation scheme and topological duality theory.
It is shown that all the standard quantum groups comprise 
holomorphic deformations. Furthermore we show that quantizing
the function algebra of a (Poisson) Lie group and quantizing 
its universal enveloping algebra are topologically dual procedures indeed.
Thus holomorphic deformation theory seems to be the appropriate
language in which to describe quantum groups as deformed Lie groups or 
Lie algebras. 
\end{abstract}
\tableofcontents
\section*{Introduction}
\addcontentsline{toc}{section}{Introduction}
In this paper we propose a new deformation scheme which we call 
{\bf holomorphic deformation} and which seems to recapture
what is actually done in the context of describing quantum 
groups as deformed Lie groups or Lie algebras. Although it
appears to be new as an explicitely formulated concept, 
we are convinced that our holomorphic deformation theory is in fact 
very close to many aspects of existing deformation procedures
in mathematical physics.

The reason to consider holomorphic deformations instead 
of the by now classical formal deformations of {\sc Gerstenhaber} 
(cf. \cite{Ger:DRA,Ger:DRAIV,GerSch:ACDT}) is twofold. 
First, one likes to obtain concrete 
deformations, i.e. deformations of the structure on a given vector 
space which are defined on this vector space and not 
only on a suitable extension.  Secondly, it is a well-known fact within 
the theory of infinite dimensional Hopf algebras that one is often forced  
to change the usual tensor product and/or the 
concept of the dual space. This change is well understood  
by introducing a (locally convex) topology on 
the Hopf algebra in question as has been shown  
in the work \cite{BonFlaGerPin:HGSQGSDRPD} of {\sc P.~Bonneau, M.~Flato, M.~Gerstenhaber} and {\sc G.~Pinczon}. 
They consider nuclear Hopf algebras $H$ and  
work with the locally convex ring $\C[[T]]$ of formal power series as a 
ring extension of $\C$ in order to formulate the concept of a deformation: 
A deformation of $H$ in their setup is a certain 
Hopf algebra structure on $H_T := \C[[T]] \hat {\otimes} H$ 
over $\C[[T]]$, where 
$\hat {\otimes}$ denotes the completed $\pi$-tensor product of locally 
convex spaces. The change of definition we propose 
(in sections \ref{HoDeLoCoAl} and \ref{HoDeNuHoAl})
is simply to replace  $\C[[T]]$ by the locally convex algebra 
${\cal O}(\Omega)$ of holomorphic functions on 
a domain $\Omega \subset \C$ (or a complex manifold $\Omega$):
A holomorphic deformation of $H$ thus is a certain Hopf algebra structure
on $H_{\Omega} := {\cal O}(\Omega) \hat {\otimes} H \cong {\cal O}(\Omega,H)$ 
over ${\cal O}(\Omega)$, where ${\cal O}(\Omega,H)$ denotes the 
locally convex space of holomorphic functions on $\Omega$ with values in $H$.
The advantage of this approach lies, among other things, in the fact that 
it is possible to actually insert values $z \in \Omega$ into the 
holomorphic deformation in order to get a deformed Hopf algebra structure 
on $H$ (and not merely on $H_T$ resp. $H_{\Omega}$). 
Furthermore, the structure maps of our concrete deformations of $H$ are 
evaluations of mappings depending holomorphically on $z$.   

Of course, in order to show that this variation of a deformation concept 
is reasonable and useful, one has to give interesting examples.
In the present paper we will show that a large part of the actually 
studied deformations of Hopf algebras, in particular those arising
in the context of quantum groups, can in fact be interpreted as being 
holomorphic deformations. This will be done in 
the final section of the present paper where we show that the Drinfeld 
and the FRT models can be regarded as being holomorphic deformations 
of the respective Lie algebras or Lie groups. 

In the course of preparing this paper we also found it useful 
to discuss the problem of what kind of locally convex topologies on the
Hopf algebra in question one should consider and whether or not there 
exist interesting and natural such topologies. 
These questions are dealt with in the second section after we have
introduced in the first section some basic definitions concerning
the compatibility of algebraic structures on a vector space with 
a given locally convex topology. 
The rather elementary considerations concerning 
locally convex topologies compatible with a given Hopf algebra, 
leads to more general topologies than the ones studied in 
\cite{BonFlaGerPin:HGSQGSDRPD}. Our discussion in section \ref{NLCNHA} 
furthermore shows that 
-- contrary to what might be expected -- introducing a topology 
on a Hopf algebra $H$ does not necessarily require an additional structure: 
There exist useful and natural topologies on $H$, even nuclear ones, 
which are completely determined by the algebraic data of $H$. 

One of the most useful topologies on a given (Hopf) algebra $H$ is, 
for instance, the locally convex projective limit topology with respect to 
the family of all finite dimensional representations of $H$. 
Working with this projective limit topology indicates the close connection of 
a holomorphic deformation with the concept of a formal deformation. 
Indeed, given a formal deformation of e.g.~a universal enveloping 
algebra ${\cal U} \g$ of a Lie algebra $\g$, it is in general not allowed to 
insert numerical values replacing the formal parameter. 
But after fixing a finite dimensional representation $\rho$ of ${\cal U} \g$ 
this can often be done with respect to the finite dimensional algebra 
$\rho({\cal U}\g)$. In particular this procedure works for deformations
considered in inverse scattering theory or conformal field theory.
Now regarding these important examples it is natural to 
introduce the projective limit topology of finite dimensional representations.
Investigating the reason why it is possible to evaluate 
the deformation parameter in the representation $\rho({\cal U} \g)$ 
at any element of a complex domain one discovers  holomorphic expressions. 
We thus arrive at a holomorphic deformation of $ {\cal U} \g$, or more 
precisely, of the completion of ${\cal U} \g$ with respect to the 
projective limit topology of all finite dimensional representations.

Our work is motivated by the desire to understand physicists work on 
deformation quantization and inverse scattering, in particular the papers
\cite{FocRos:PSMFCRSRM} by {\sc Fock, Rosly} and 
\cite{AleGroSch:CQHCSTI,AleGroSch:CQHCSTII} by 
{\sc Alekseev, Grosse, Schomerus},
where algebras are ``concretely deformed'' by a real parameter $\hbar$ or 
$q$ and not only by a formal one. 
We hope that the concept of holomorphic deformations will give the means
to better understand the deformation quantization of the moduli 
space of flat connections on a given connected, oriented and compact surface 
with marked points as described in
\cite{FocRos:PSMFCRSRM,AleGroSch:CQHCSTI,AleGroSch:CQHCSTII},
and to compare it with the quantization in 
{\sc Scheinost, Schottenloher} \cite{ScheScho:MQMSFPB}.

Let us also mention that our notion of a deformation creates the means
to consider formal and holomorphic deformations in one common language. 
For more details about this approach we refer the interested reader to 
{\sc Pflaum} \cite{Pfl:LADQ,Pfl:NCDQNOQCB}. 

The main ingredient to our theory of holomorphic deformations is not so 
much an algebraic viewpoint, but a functional analytic one. In particular
we will often use the notions of locally convex spaces, nuclearity, and 
topological tensor products. For the convenience of the reader we 
therefore explain in the appendix the most important functional analytic 
concepts used in this article. 

The problem of classifying holomorphic deformations is not considered 
here. This will be the subject of a forthcoming paper together with the 
study of a suitable cohomology theory.
\\[3mm]
{\bf Acknowledgement.} An essential part of this work has been done while
both of the authors participated for two weeks in the {\it Research in Pairs}
program at the Forschungsinstitut Oberwolfach. We thank the staff at 
Oberwolfach for the excellent working conditions and greatly
acknowledge financial support by the Stiftung Volkswagenwerk.

\section{Topological and algebraic structures}
\label{NuSpAlSt}

In order to formulate our deformation concept in the next section
we need a topological structure on a given vector space. Therefore, in this
section we study locally convex structures on a vector space $E$ which 
are compatible with additional algebraic structures on $E$ such as the structure
of an algebra, of a coalgebra or of a Hopf algebra. A collection of 
definitions and results on locally convex spaces is provided in the 
Appendix.

In the following $\K$ will denote one of the fields $\R$ and $\C$
together with the Euclidean topology. Furthermore, the locally 
convex topologies
on a given $\K$-vector space $E$ are always assumed to be Hausdorff 
and complete.  The collection of the (complete and Hausdorff) 
locally convex spaces together with the
$\K$-linear continuous maps form a category $\cat{Lcs}$.
For two locally convex spaces $E$ and $F$
the completion  of the tensor product 
$E \otimes F$ endowed with the $\pi$-topology is denoted 
by $E \hat{\otimes} F$. The significance of $E \hat{\otimes} F$ is 
apparent by the following universal property:
There is a natural bijection between the 
space of continuous bilinear maps $E \times F \rightarrow G$ into a third
locally convex space $G$ and the space of continuous linear maps 
$E \hat{\otimes} F \rightarrow G$. With $\hat{\otimes}$ as tensor product 
functor $\cat{Lcs}$ becomes a symmetric monoidal category.

Furthermore, a nuclear space is always supposed to be a locally 
convex $\K$-vector space $E$ which is Hausdorff, complete  and nuclear. 
The completion $E \hat{\otimes} F$ of $E \otimes F$ for two 
nuclear spaces $E$ and $F$ is again a nuclear space
(see the Appendix for precise definitions and examples as well as
{\sc Pietsch} \cite{Pie:NLCS} and {\sc Tr\`eves} \cite{Tre:TVSDK}
for proofs and details about nuclear spaces). Thus nuclear spaces 
form a monoidal subcategory $\cat{Nuc}$ of $\cat{Lcs}$. 

We call a nuclear space $E$ {\bf strictly nuclear} 
if its (strong) dual $E'$ is nuclear as 
well, if $E$ is reflexive (i.e. the strong dual of $E'$ is 
isomorphic to $E$ as a locally convex space) and if it fulfills the 
{\bf duality condition} (cf. (\ref{IsStNuSp}) in the Appendix), 
i.e.~the canonical linear 
mapping $E' \otimes E' \rightarrow (E \hattimes E)'$
extends to an algebraic and topological isomorphism 
$$E' \hattimes E' \rightarrow (E \hattimes E)'.$$
Nuclear Fr\'echet spaces or duals of nuclear 
Fr\'echet spaces are strictly nuclear
(cf.~{\sc Tr\`eves} \cite{Tre:TVSDK}) as well as nuclear LF-spaces.
Restricting the object class of $\cat{Nuc}$ to strictly nuclear spaces
we receive a proper subcategory $\cat{sNuc}$ of $\cat{Nuc}$. 

Algebraic structures can now be formulated within $\cat{Lcs}$, (resp.  
$\cat{Nuc}$ or $\cat{sNuc}$) in order to 
obtain locally convex (resp. nuclear or strictly nuclear) 
algebras,  Hopf algebras, etc.
In this section we will give detailed definitions for these objects and 
will examine them. In the next section we present natural examples.

\begin{Definition}

A {\bf locally convex algebra} is a locally convex  space $A$ together 
with continuous linear mappings $\mu :\:A \hat{\otimes} A \rightarrow A$ 
and $\eta:\:\K \rightarrow A$
such that $\mu$ fulfills the associativity constraint $\mu \circ (\mu 
\hattimes \id_A ) = \mu \circ (\id_A \hattimes \mu)  $ and $\eta$ gives
rise to a unit: 
$\mu \circ (\id_A \hattimes \eta) = \mu \circ (\eta \hattimes \id_A) \cong 
\id_A $. A homomorphism between locally convex algebras $A$ and
$\tilde{A}$ is just a continuous linear map $f: A \rightarrow \tilde{A}$
such that $\tilde{\mu} \circ ( f \hattimes f) = f \circ \mu$ and
$\tilde{\eta} \circ f = \eta$. 

A {\bf locally m-convex algebra} is a locally convex algebra $A$ for 
which there exists a defining family of multiplicative seminorms (see 
Appendix).

A {\bf nuclear algebra} (resp. a {\bf strictly nuclear algebra} ) 
is a locally convex algebra for which the underlying
locally convex space $A$ is a nuclear space (resp. a strictly nuclear
space).

Similarly one defines the concepts of {\bf locally convex coalgebra}, 
{\bf locally convex  bialgebra} and  {\bf locally convex Hopf algebra}: 
These are locally convex spaces together with appropriate 
continuous structure
maps such that, respectively, the axioms of coassociativity, counit,
compatibility of multiplication and comultipication and the axiom of the
antipode are fulfilled with $\hattimes$ as tensor product functor.
Likewise there exists an appropriate notion of morphism of locally 
convex coalgebras, locally convex bialgebras and locally convex 
Hopf algebras, i.e. continuous linear maps which leave the 
structure maps invariant.

Finally, a ({\bf strictly}) {\bf nuclear Hopf algebra} $H$ is a locally convex 
Hopf algebra such that the underlying locally convex space $H$ is (strictly)
nuclear. Similarly define the notions of
({\bf strictly}) {\bf nuclear coalgebras} and {\bf bialgebras}.
\end{Definition}
\begin{rRemark}
  A nuclear space $A$ comprises a nuclear algebra if and only if 
  it has an underlying structure of a $\K$-algebra such that  
  the multiplication $\mu : A \times A \rightarrow A$ is
  continuous. 
  A similar result does not hold for coalgebras. Namely
  there exist nuclear coalgebras $C$ which do not have an 
  underlying structure of a coalgebra. In other words this means that the
  map $\Delta : C \rightarrow C \hat{\otimes} C$ need not have
  its image in $C \otimes C$. An example is given by the quantized
  ${\fr sl} (N+1, \C)$ of section \ref{QuGrHoHoAl}.
\end{rRemark}
Within the topological setting we also have to define appropriate
topological versions of R-matrix, triangularity, coquasitriangularity, etc.
\begin{Definition}
  Let $H$ be a locally convex Hopf algebra or bialgebra. It is called 
  {\bf topologically quasitriangular} if there exists an invertible element
  ${\cal R} \in H \hattimes H$ such that the following conditions hold:
\begin{eqnarray}
  \tau \circ \Delta (a) & = & {\cal R} \, \Delta (a) \, {\cal R}^{-1}
  \\[2mm]
  (\Delta \otimes {\rm id}) ({\cal R} ) & = & {\cal R}_{13} \, {\cal R}_{23}
  \\[2mm]
  ({\rm id} \otimes \Delta) ({\cal R} ) & = & {\cal R}_{13} \, {\cal R}_{12},
\end{eqnarray}
where $\tau : H \hattimes H \rightarrow H \hattimes H$ is the flip-morphism
and ${\cal R}_{12},\, {\cal R}_{13},\, {\cal R}_{23}$ are the obvious
extensions of ${\cal R}$ to $H \hattimes H \hattimes H$ which are trivial
on the third, resp.~second, resp.~first factor. In this case ${\cal R}$
is called the {\bf topological universal R-matrix} of $H$.
If additionally 
\begin{eqnarray}
 {\cal R}^{-1} & = & \tau \circ {\cal R}
\end{eqnarray}
then $H$ is called {\bf topologically triangular}.

Dually $H$ is called {\bf topologically coquasitriangular}, if
there exists a continuous bilinear map $< \: | \: > : \:  H \hattimes H 
\rightarrow \C$,
the {\bf braiding form}, such that for all $a,b,c \in H$
\begin{eqnarray}
  \sum \limits_{(a),(b)} \, < a_1 | b_1 > \, a_2 \, b_2 & = & 
  \sum \limits_{(a),(b)} a_1 \, b_1 \, < a_2 | b_2 >, \\[2mm]
  < a | bc > & = &  \sum \limits_{(a)} <a_1 | b> \, < a_2 | c> , \\[2mm]
  < ab | c> & = &  \sum \limits_{(c)} <a | c_1> \, < b | c_2>.
\end{eqnarray}

\end{Definition}

One of the main reasons to consider strictly nuclear Hopf algebras instead of
just Hopf algebras lies in the fact that the category of strictly
nuclear Hopf algebras has duals.
\begin{Proposition}
\label{PrDu}
  Let $H$ be a strictly nuclear Hopf algebra. Then the dual space $H'$ carries
  in a natural way the structure of a nuclear Hopf algebra.
  Moreover, $H$ can be recovered as $H''$.

  If $H$ is topologically quasitriangular (resp.~coquasitriangular) then
  $H'$ is topologically coquasitriangular (resp.~quasitriangular).

  The same properties hold for reflexive locally convex Hopf algebras 
  fulfilling the duality condition.
\end{Proposition}
\begin{Proof}
  For the first part of the proposition just apply the isomorphism 
  $(H \hat{\otimes} H)' \cong H' \hat{\otimes} H'$
  to obtain the coproduct $\Delta'$ on $H'$ as the pull-back 
\begin{equation}
  \mu^* : \: H' \longrightarrow (H \hat{\otimes} H)', 
  \quad f \longmapsto f \circ \mu.
\end{equation} The other structure maps of $H'$ are directly defined by    
  transposition.

  For the second part note that an element ${\cal R} \in H \hattimes H$ 
  induces a continuous bilinear form $< \: | \: >_{\cal R}: \,
  H ' \hattimes H ' \rightarrow \C$ by
\begin{equation}
  f \otimes g \mapsto f \otimes g ({\cal R}).
\end{equation}   Furthermore, by the isomorphism 
  $(H \hat{\otimes} H)' \cong H' \hattimes H' $
  every continuous bilinear form ${<\: | \: >:} \, H \otimes H \rightarrow \C$
  can be interpreted as an element ${\cal R}_{<\, |\, >} \in H' \hattimes H'$.  
  The proof of the required algebraic properties for the thus defined 
  braiding form $< \: | \: >_{\cal R}$ resp.~R-matrix 
  ${\cal R}_{<\, |\, >}$ follows exactly like in the well-known 
  finite dimensional case. 
\end{Proof}
Let us finally mention that all the nuclear spaces given in the 
Appendix, in particular the function spaces ${\cal E} (\Omega)$ and 
${\cal O} (\Omega )$ as well as all finite dimensional algebras 
(resp.~coalgebras, bialgebras and Hopf algebras), comprise
examples of nuclear algebras (resp.~coalgebras, bialgebras and Hopf algebras).

\section{Natural locally convex and nuclear Hopf algebras}
\label{NLCNHA}
In this section we will consider some general constructions and 
examples of locally convex Hopf algebras. 
Note that analoguous topological constructions can be carried through
for locally convex algebras and bialgebras as well.
\begin{rPar} {\bf Inductive limit topologies on Hopf algebras.}
\label{InLiToHoAl}
  On a given Hopf algebra $H$ over $\K$ we  
  can always consider the finest locally convex topology. 
  Then $H$ is a complete locally convex Hausdorff space and, since all
  structure maps are automatically continuous for the finest locally 
  convex topology, $H$ is a locally convex Hopf algebra. But this locally
  convex Hopf algebra $H$ is locally m-convex only if $H$ is finite 
  dimensional and nuclear only if $H$ is of countable dimension. 
  The dual $H'$ with the strong topology carries the
  coarsest locally convex topology. Since $H \cong \K^{(\Lambda)}$ 
  satisfies the duality  condition
  (see Appendix (\ref{IsKL})) $H'$  is, according to 
  proposition \ref{PrDu}, a locally convex Hopf algebra as well.
  For the special case of a group algebra
  $H = \K G$ of a group $G$ this dual is the Hopf algebra $H' \cong \K^G$ of
  all functions on $G$.
\end{rPar}
\begin{rPar}{\bf Projective limit topologies on Hopf algebras.}
\label{PrLiToHoAl}
\label{ToSeSiUnEnAl}
  Alternatively  one could provide a given 
  Hopf algebra $H$ over $\K$ with the locally convex 
  projective limit {\bf topology of finite dimensional representations}, i.e.
  with the coarsest locally convex topology leaving continuous all 
  finite dimensional representations 
  $\varphi: H \rightarrow \End V$. 

    Let us assume that  these representations   
  separate the points of $H$. 
  According to the later proved Proposition
  \ref{PrSeFiRe} this is e.g. the case for finitely generated Hopf algebras 
  $H$. It also holds for the universal enveloping algebra $H = {\cal U} \g$ 
  of a finite dimensional Lie algebra $\g$. Then $H$ is Hausdorff locally 
  convex space. $H$ is complete only if $H$ is finite dimensional. The
  completion $\hat{H}$ of $H$ however lies in $\cat{Nuc}$ since it is 
  a locally convex projective limit of finite dimensional spaces.
  All structure maps of the Hopf algebra $H$ are continuous since the topology
  is adapted to the finite dimensional representations. Therefore, they can 
  be uniquely extended to $\hat{H}$ and thus turn $\hat{H}$ into a nuclear
  Hopf algebra. In addition $\hat{H}$ is locally m-convex.

    The dual $\hat{H}'$  of $\hat{H}$ (or of $H$) is the space of 
  matrix coefficients on $H$:
\begin{eqnarray}
  \hat{H}' & = & \left\{ \xi \circ \hat{\varphi}  \, | \: 
  \varphi : H \rightarrow
  \End V \mbox{ finite dim.~representation, }
  \xi \in ({\rm End} \, V)' \right\} .
\end{eqnarray}
  Hence, $\hat{H}'$ coincides with the restricted dual $H^{\circ}$ of $H$
  (see {\sc Chari, Pressley} \cite{ChaPre:GQG} Chapter 4.1.D).
  The strong topology on $\hat{H}'$ is given by the locally convex inductive 
  limit topology of the maps
\begin{equation}
  \varphi':
  ({\rm End} \, V)' \rightarrow \hat{H}', \quad 
  \xi \mapsto \xi \circ \hat{\varphi},
\end{equation}
  where $\varphi$ runs through all finite dimensional representations of
  $H$. Thus the strong topology on $\hat{H}'$ is the finest locally convex 
  topology. Although the locally convex space $\hat{H}' = H^{\circ}$ 
  is in general not nuclear, it satisfies the duality condition 
  (cf. Appendix, Eq.~(\ref{IsKL})). Therefore, as in the Proposition \ref{PrDu}
  the transpositions of the 
  structure maps of the locally convex Hopf algebra $\hat{H}$ define the  
  structure of a locally convex Hopf algebra on $\hat{H}'$. Dualizing again
  one gets the Hausdorff completion ${H^{\circ}}'$ of $H$ (which is 
  $\hat{H}$ if $H$ is assumed to be Hausdorff).

    The locally convex space $\hat{H}'$ will 
  be nuclear, even strictly nuclear, whenever 
  countably many of the finite dimensional representations generate 
  the topology of $H$, i.e. if $\hat{H}$ is Fr\'echet. This is the case
  e.g. for the universal enveloping algebra ${\cal U} \g$ of a 
  semi-simple Lie algebra $\g$ (see below).

    In the same spirit one can attach to $H$ the projective limit topology
  with respect to all homomorphisms $H \rightarrow A$, where
  $A$ is a nuclear locally m-convex Fr\'echet algebra. We call the resulting 
  projective limit topology the 
  {\bf topology of nuclear Fr\'echet representations}.
  The completion 
  $\check{H}$ of $H$ with respect to this topology again is a nuclear Hopf 
  algebra. Furthermore, one has a natural continuous inclusion 
  $\check{H} \rightarrow \hat{H}$. 

    More generally, one can consider other projective systems of representations
  of $H$ in order to define appropriate locally convex topologies on $H$ which 
  have their origin in purely algebraic properties of $H$.
\begin{rRemark}
\label{C[T]}
  In the case $A = \C [ x_1,...,x_n]$ the topology of finite dimensional
  representations is the topology of pointwise convergence
  of all derivatives. Hence, $\hat{A}$ is not a Fr\'echet space.
  Similarly, the tensor algebra ${\rm T} V$ endowed with the topology of 
  finite dimensional representations is a nuclear m-convex algebra which 
  is not metrizable.

  However, the completion of $A= \C [ x_1,...,x_n]$ with respect to 
  the topology of nuclear Fr\'echet representations is isomorphic to 
  the nuclear m-convex Fr\'echet algebra
  ${\cal O } (\C^n)$ of entire holomorphic functions on $\C^n$.
\end{rRemark}

  We will determine in the following  the topologies of finite dimensional
  representations and of Fr\'echet representations for some important 
  algebras.
\begin{Proposition}
\label{PrSeFiRe}
  Let $A$ be a finitely generated algebra and ${\rm T} V 
  \stackrel{\pi}{\rightarrow} A$
  a presentation of $A$, where ${\rm T} V$ is the tensor algebra of a 
  finite dimensional $\K$-vector space $V$. 
  Denote by $\hat{\rm T} V$ and $\hat{A}$ (resp.~$\check{\rm T} V$ and 
  $\check{A}$) the completions of ${\rm T} V$ and $A$ with repsect to the 
  topology of finite dimensional representations (resp.~of nuclear Fr\'echet 
  representations).
  Then the algebras $\hat{\rm T} V$, $\hat{A}$, $\check{\rm T} V$ and 
  $\check{A}$ are locally m-convex nuclear Hausdorff spaces. The spaces
  $\check{\rm T} V$ and $\check{A}$ are even Fr\'echet.
  Furthermore, the presentation ${\rm T} V \stackrel{\pi}{\rightarrow} A$
  extends uniquely to surjective and open maps 
  $\hat{\rm T} V \stackrel{\hat{\pi}}{\rightarrow} \hat{A}$ and
  $\check{\rm T} V \stackrel{\check{\pi}}{\rightarrow} \check{A}$,
  i.e.~to topological presentations of $\hat{A}$ and $\check{A}$.
\end{Proposition}
\begin{Proof}
  By the universal property of the complete hull we have unique
  homomorphisms $\hat{\rm T} V \stackrel{\hat{\pi}}{\rightarrow} \hat{A}$ and
  $\check{\rm T} V \stackrel{\check{\pi}}{\rightarrow} \check{A}$
  both extending ${\rm T} V \stackrel{\pi}{\rightarrow} A$.
  We will show that they have the claimed properties.
  
  First consider the topology of finite dimensional representations.
  Let us show that this topology is Hausdorff or in other words that the
  finite dimensional representations of $A$ separate the points of $A$. 
  Consider the ideals $I_n = \oplus_{k \geq n}
  V^{\otimes k}$ in ${\rm T} V$. Their images in $A$ define ideals $J_n$
  in $A$. 
  Now for two elements $a,b \in A$, $ a\neq b$ there exists an $n \in \N$ 
  large enough such that $a - b $ dos not vanish in $A / J_n$. 
  But $A / J_n$ is finite dimensional and an $A$-module. 
  Thus the points $a$ and $b$ are separated
  by the representation of $A$ on $A / J_n$.
  
  The continuous homomorphism $\hat{\rm T} V \rightarrow \hat{A}$ is an open
  map. To see this choose an open set $U \subset \hat{\rm T} V$. We can assume
  that there exists a finite dimensional representation 
  $\varphi : \hat{\rm T} V \rightarrow \End W$ and an open 
  $O \subset \End W$ such that $U = \varphi^{-1} (O)$. Let $\hat{I}$ be the 
  kernel of $\hat{\pi}$, and $\tilde{W}$ be the algebra 
  $\im{\varphi} / \varphi \, (\hat{I})$. 
  Then $\varphi$ induces a representation
  $\tilde{\varphi} : A \rightarrow \tilde{W} \subset \End \tilde{W}$.
  As projections between finite dimensional spaces are open there exists
  an open $\tilde{O} \subset \End \tilde{W}$ such that 
  $\tilde{O} \cap \im{\tilde{\varphi}} = O + \varphi (\hat{I})$.
  We then have $\hat{\pi} (U) = \tilde{\varphi}^{-1} (\tilde{O})$, 
  i.e.~$\hat{\pi} (U)$ is open in $\hat{A}$. 
  
  As $\hat{\rm T} V \rightarrow \hat{A}$ is continuous, open and has
  dense image, it is surjective, hence, provides a continuous presentation
  of $\hat{A}$.  \vspace{2mm} 

  Now let us consider the case of nuclear Fr\'echet representations.
  Choose a basis $(x_1,...,x_n)$ of $V$ and recall that ${\rm T} V$
  is canonically isomorphic to the linear space 
  $\C^{(<x_1,...,x_n>)}$, where  ${<x_1,...,x_n>}$ is the free half group
  generated by $x_1,...,x_n$. For any $n$-tupel $\alpha \in (\Q^+)^n$
  define the seminorm $p_{\alpha}$ on $V$ by
\begin{equation}
  p_{\alpha} \left( \sum \limits_{1 \leq k \leq n} \, \lambda_k \, x_k \right)
  \: = \:  \sum \limits_{1 \leq k \leq n} \, | \lambda_k | \, \alpha_k.
\end{equation}
  This seminorm naturally extends to a seminorm on $V^{\otimes k}$, $k \in \N$
  and then to a seminorm on ${\rm T} V$. By definition 
  $ p_{\alpha} ( v \otimes w ) = p_{\alpha} ( v) \,   p_{\alpha} ( w )$
  holds for every $v,w \in {\rm T} V$. 
  Hence, the  countable system of seminorms 
  $p_{\alpha}$ gives rise to a metrizable locally m-convex topology on 
  ${\rm T} V$. Consequently the completion $E$ of ${\rm T} V$ with respect 
  to the topology generated by the seminorms $p_{\alpha}$ is locally m-convex
  and  Fr\'echet. Now assume that $\rho: \, {\rm T} V \rightarrow F$ 
  is a representation with $F$ a nuclear locally
  m-convex Fr\'echet algebra. Let $p$ be a multiplicative seminorm on $F$.
  Choose $\alpha_k \in \Q^+$, $k =1,...,n$ so large that 
  $ p (\rho (x_k)) \leq p_{\alpha} ( x_k) $. By the multiplicativity
  of $p$ and the definition of $p_{\alpha}$ this implies 
  $p (\rho (v )) \leq p_{\alpha} (v)$ for every $v \in {\rm T} V$.
  Hence, $\rho$ extends to a continuous homomorphism $E \rightarrow F $.
  
  Let us suppose for a moment that $E$ is nuclear.
  Since the inclusion ${\rm T}V \rightarrow E$ is a continuous Fr\'echet 
  representation, the nuclearity of $E$ entails the relation 
  $E \cong \check{T} V$. Denoting by $\check{I}$ the completion
  of the kernel $I$ of $\pi$, the quotient algebra $E / \check{I}$ then
  is nuclear as well, locally m-convex and Fr\'echet. 
  As $A$ lies densely in $E / \check{I}$ the same argument like for $E$ shows 
  that $E / \check{I}$ is the completion $\check{A}$. Hence, the  
  claim follows.
 
  So, it remains to show that $E$ is nuclear. We will achieve this
  by an argument using the sequence spaces of {\sc Pietsch} \cite{Pie:NLCS}. 
  Denote for every $\beta \in <x_1,...,x_n>$ by $\ell (\beta) $ the 
  multilength of $\beta$; that means $\ell (\beta)$ is the $n$-tupel
  $(\ell_1,...,\ell_n)$ where $\ell_k$ counts how many times $x_k$ appears 
  in $\beta$. The absolute length $\ell_1 + ... + \ell_n$ of $\beta$
  is denoted by $|\ell (\beta ) |$. 
  Now choose a bijection $\xi: \N \rightarrow  <x_1,...,x_n>$ 
  such that for $\beta, \tilde{\beta} \in <x_1,...,x_n>$ fulfilling 
  $| \ell (\beta) | < | \ell (\tilde{\beta}) | $ the relation 
  $\xi^{-1} (\beta) < \xi^{-1} (\tilde{\beta})$ holds.
  By definition $E$ is then isomorphic to the sequence space 
  (cf.~{\sc Pietsch} \cite{Pie:NLCS}, Chapter 6) 
  generated by the
  system of sequences $(\lambda_k^{\alpha} )_{k \in \N}$ with 
  $\lambda_k^{\alpha} = \alpha^{\ell (\xi (k))}$ and
  $\alpha \in (\Q^+)^n$. For $ 0 < q < \frac{1}{n}$ we now have
\begin{equation}
 \sum \limits_{\beta \in <x_1,...,x_n> } \, q^{| \ell (\beta ) |}
 \: = \: \sum \limits_{k \in \N} \, ( n \, q) ^k \: < \: \infty.  
\end{equation}     
  By {\sc Pietsch} \cite{Pie:NLCS}, Theorem 6.1.2. this implies that $E$ is nuclear.
  This proves the claim. 
\end{Proof}
\end{rPar}
\begin{rPar} {\bf Matrix coefficients of group representations.}
\label{MaCoGrRel}
  For a group $G$ let us consider the Hopf algebra of 
  complex-valued matrix coefficients
  ${\cal R}_0 (G) := \left\{ \xi \circ \varphi |\:
  \varphi : G \rightarrow \GL V \mbox{ finite dimensional representation, }
  \xi \in ({\rm End} \, V)' \right\}$ in the light of \ref{PrLiToHoAl}
  and \ref{InLiToHoAl}. (Here we assume $V$ always to be finite 
  dimensional {\bf complex} vector spaces.) 
  ${\cal R}_0 (G)$ with the finest locally convex topology is a locally convex  
  Hopf algebra fulfilling the duality condition. This topology can also be 
  described as the locally convex inductive limit of the maps 
  $\varphi' : (\End V)' \rightarrow {\cal R}_0 (G) $
  where $\varphi$ runs through all finite dimensional representations of $G$. 
  
    The dual ${\cal R}_0 (G)'$ of ${\cal R}_0 (G)$ endowed with the strong  
  topology is a locally convex projective limit 
  of finite dimensional algebras and thus is a nuclear Hopf algebra which is
  locally m-convex. By the duality condition for ${\cal R}_0 (G)$ the
  dual ${\cal R}_0 (G)'$ obtains as in Proposition \ref{PrDu} the structure 
  of a nuclear Hopf algebra.
 
    In case of a topological group $G$ we replace ${\cal R}_0 (G)$ by the
  continuous matrix coefficients ${\cal R} (G) \subset {\cal R}_0 (G)$.
  ${\cal R} (G)$ with the finest locally convex topology is a locally convex   
  Hopf algebra as well and the dual ${\cal R} (G)'$ is a 
  nuclear m-convex algebra.
  In general,  ${\cal R} (G)$ is not nuclear. For compact groups, however,
  ${\cal R} (G)$ is strictly nuclear since by the theorem of Peter and Weyl 
  it is of countable dimension. Moreover, the dual ${\cal R} (G)'$ is a 
  Fr\'echet nuclear locally m-convex algebra.

  The dual ${\cal R}_0 (G)'$ contains all the evaluations 
  $\delta_x : {\cal R}_0 (G) \rightarrow \C$, $ f \mapsto f(x)$,
  where $x \in G$. The map $\delta: G \rightarrow {\cal R}_0 (G)'$ defines
  a Hopf algebra morphism $\delta :\C G \rightarrow {\cal R}_0 (G)'$. In the
  case of a topological group one  analogously has a 
  natural continuous Hopf algebra map 
  $\delta : \C G \rightarrow {\cal R} (G)'$. 
  For compact groups $\delta$ is injective, open onto its image, and 
  ${\rm span}(\delta(G))$ is dense in ${\cal R} (G)'$ (cf.  
  \cite{BonFlaGerPin:HGSQGSDRPD}).
 
  Other models of locally convex Hopf algebras consisting of representative 
  functions can be considered by fixing a suitable class of representations
  of $G$ in locally convex vector spaces being closed under dualizing, 
  (finite) direct sums and (finite) completed tensor products.
\end{rPar}
\begin{rPar} {\bf Universal enveloping algebras.}
\label{UnEnAl}
  Starting with a Lie algebra $\g$ over $\K$ the universal enveloping
  algebra \, ${\cal U} \g$ is a Hopf algebra over $\K$.
  One can study \, ${\cal U} \g$ as a topological Hopf algebra with respect
  to the finest locally convex topology as in \ref{InLiToHoAl}. 
  However, the projective limit topologies of \ref{PrLiToHoAl}
  are more interesting  in some aspects.

    So consider the coarsest locally convex topology on 
    ${\cal U} \g$  such that all finite dimensional representations 
    $\rho: {\cal U} \g \rightarrow {\rm End}_{\C} \, V$ 
    are continuous. Let us assume that these finite dimensional representations 
    separate the points of ${\cal U} \g$. This is always the case for a
    finite dimensional Lie algebra by Ado's theorem. Then the 
    topology of finite dimensional representations 
    on ${\cal U} \g$ is Hausdorff, and the completion 
    $\hat{\cal U} \g$ of ${\cal U} \g$ is a nuclear Hopf algebra as 
    in \ref{PrLiToHoAl}. Moreover, $\hat{\cal U} \g$ is locally m-convex.

    For a Lie group $G$ with Lie algebra $\g$ there is a natural
    map $i$ relating $\hat{\cal U} \g$ and the nuclear 
    (cf. \ref{MaCoGrRel})  Hopf algebra 
    ${\cal R} (G)'$. $i: {\cal U} \g \rightarrow {\cal R} (G)'$ is defined
    by $i(X) (f) = L_X f (e)$, $f \in {\cal R} (G)$, where $L_X$ is the left
    invariant differential operator on $G$ given by $X \in {\cal U} \g$ and 
    where $e$ is the unit of $G$. Comparing the topologies one sees that 
    $i$ is continuous since every finite dimensional representation 
    $\varphi$ of $G$ induces a representation $\dot{\varphi}$ on $\g$ by
    differentiation and since 
    $i(X) (\xi \circ \varphi) = \xi (\dot{\varphi}(X))$. 
    Hence, $i$ can be extended to a continuous $\R$-linear map 
    $i : \hat{\cal U} \g \rightarrow {\cal R} (G)'$ which 
    is a morphism of nuclear Hopf algebras.
    $i$ is not injective in general. For connected and simply connected 
    Lie groups, however, $i$ is injective. Moreover, in that case 
    $i$ is an open map
    onto its image $ i( \hat{\cal U} \g) $, 
    since the finite dimensional complex representations of $\g$ and $G$ 
    are in one-to-one correspondence. 
    Therefore, $\hat{\cal U} \g$   can be considered as a closed nuclear 
    sub-Hopf algebra of ${\cal R} (G)'$.  
    But in general $i$ does not need to be open for every Lie group $G$.
    Take for example $G = {\rm U}(1)$ and its Lie algebra $\g \cong \R$.
    Then $\g$ is also the Lie algebra of $\R$ and 
    ${\cal U} \g \cong \R[T]$. Though $\hat{\cal U} \g$ is a 
    closed sub-Hopf algebra of ${\cal R} (\R)' \cong {\C}^{\R}$
    the map  
    $i: {\cal U} \g \rightarrow {\cal R} ({\rm U}(1))' \cong {\C}^{\N}$
    cannot be open since the topology of finite dimensional representations 
    on ${\cal U} \g$ is not metrizable (cf. \ref{C[T]}).
\end{rPar}

\begin{rPar} {\bf Simple Lie algebras.}
\label{SiLiAl}
    If $\g$ is a simple Lie algebra over $\C$ then finitely many
    of the finite dimensional representations already generate all finite    
    dimensional representations of $\g$ (via finite sums and tensor 
    products; e.g. for $\g = {\fr sl}(N,\C)$ the finite 
    dimensional representations are generated
    by the fundamental representation 
    ${\fr sl} (N,\C) \subset {\fr gl} (N,\C)$). 
    As a consequence, the topology of finite dimensional representations on 
    $\hat{\cal U} \g$ is metrizable and hence Fr\'echet. Therefore, in the
    simple case $\hat{\cal U} \g$ is in particular strictly nuclear.
    
    For a simple complex Lie algebra $\g$ with corresponding connected
    and simply connected Lie group $G$ the image 
    $ i({\cal U} \g) $ is dense in ${\cal R} (G)'$. Hence, by the above, 
    the map $i : \hat{\cal U} \g \rightarrow {\cal R} (G)'$
    is an isomorphism of Fr\'echet algebras. As a consequence, in this
    situation we have a natural complete duality between 
    $\hat{\cal U} \g$ and  ${\cal R} (G)'$. In particular, 
    $\hat{\cal U} {\fr sl} (N, \C) \cong {\cal R} ({\rm SL} (N,\C)'$.

    In the case of a compact Lie group $G$ the map 
    $i: {\cal U} \g \rightarrow {\cal R} (G)'$  is injective as well.
    It can be continued to a $\C$-linear injective map 
    $i: {\cal U} \g^{\C} \rightarrow {\cal R} (G)'$ by complexification
    of the Lie algebra $\g$. 
    The image $ i \left( {\cal U} \g^{\C}\right) $ turns out 
    to be dense in  ${\cal R} (G)'$. However, the induced topology on  
    $ i \left( {\cal U} \g^{\C} \right) \subset {\cal R} (G)'$ 
    does in general not
    coincide with the topology coming from the projective topology of finite
    dimensional representations; see, e.g. the example of $U(1)$. Instead
    of this, the inclusion $i$ induces a new locally convex topology 
    on ${\cal U} \g^{\C}$ which can be described as the locally convex
    projective limit of all representations $\dot {\varphi}$ which are 
    derivatives of finite dimensional continuous representations 
    $\varphi$ of the 
    group $G$. This topology depends on the group in question and not 
    only on the Lie algebra $\g$. It is always metrizable and nuclear. 
    Hence, the completion --  which we denote by $\tilde{\cal U} \g^{\C}$
    -- is a Fr\'echet nuclear Hopf algebra naturally isomorphic to 
    ${\cal R} (G)'$.

    As the completion $\check{\cal U} \g $ of ${\cal U} \g $ with respect
    to the topology of nuclear Fr\'echet representations  naturally lies in
    $\hat{\cal U} \g$ we have a continuous inclusion 
    $\check{\cal U} \g  \rightarrow {\cal R} (G)'$ as well. 
    Note that $\hat{\cal U} \g$,  $\check{\cal U} \g $,  
    $\tilde{\cal U} \g^{\C}$ and ${\cal R} (G)'$ are locally m-convex 
    algebras as projective limits of locally m-convex algebras.
\end{rPar}

%
\begin{rPar} {\bf Function and distribution spaces on Lie groups}
\label{FuDiSpLiGr}
  All the above examples of locally convex and nuclear Hopf algebras 
  are based essentially on the 
  algebraic structures of the Hopf algebra in question. In many cases they 
  carry the finest or the coarsest locally convex topology. In the case of a
  Lie group $G$ other topological structures can be imposed on certain
  Hopf algebras of smooth functions on $G$. 
  These topologies arise from the analytical structure on $G$. 
  Similarly new locally convex Hopf algebra structures
  on the universal enveloping algebra ${\cal U} \g$ of the Lie algebra $\g$ 
  are of interest.

  For example, on the space of complex valued smooth functions 
  ${\cal E} (G) := {\cal C}^{\infty} (G)$ the topology of uniform convergence 
  of all derivatives on the compact sets of $G$ induces on ${\cal E} (G)$ 
  the structure of a strictly nuclear Hopf algebra. ${\cal E} (G)$ 
  is a locally m-convex
  algebra. Its dual ${\cal E}'(G)$ contains $G$ by the evaluations 
  $\delta_x$, $x \in G$, and hence it contains the group algebra 
  $\C G$. ${\cal E}' (G)$ also contains ${\cal U} \g$ via the injective
  Hopf algebra map $i : {\cal U} \g \rightarrow
  {\cal E}' (G)$, where $i(X) = \delta_e \circ L_X$ now acts on the bigger 
  space ${\cal E} (G)$. Note that this map is not continuous in general. It
  gives rise to yet another locally convex topology on ${\cal U} \g $ whose
  completion is a strictly nuclear Hopf algebra.

  Another dual pair of strictly nuclear Hopf algebras is given 
  by the test functions
  ${\cal D} (G) = {\cal C}_{c}^{\infty} (G) \subset {\cal E} (G)$ and its
  dual ${\cal D}'(G)$, the space of distributions on $G$.
  We now have the following continuous inclusions:
\begin{eqnarray}
\begin{array}{ccccc}
  {\cal R} (G) & \subset & {\cal E} (G) \,\, , \,\,
  {\cal E}' (G) & \subset & {\cal R} (G)' \\
  \K G & \subset & {\cal E}' (G) & \subset & {\cal D} ' (G) \\
  {\cal D} (G) & \subset & {\cal E} (G) & \subset & \K^G \cong \K G'. 
\end{array}
\end{eqnarray}
 
\end{rPar}
\begin{rPar}{\bf FRT-bialgebras.}
\label{FRT}
 In the following we will show how the construction of quantum matrix algebras
 by {\sc Faddeev, Reshetikhin, Takhtajan} \cite{FadResTak:QLGLA} 
 can be implemented in our concept
 of nuclear Hopf algebras and comprises a further nontrivial example. 

 First we introduce some notation. Let $A$ be an arbitrary algebra. The
 space of $n \times m$ matrices with entries from $A$ is denoted
 by ${\rm M} (n \times m, A)$. The Kronecker product is then the mapping
\begin{equation}
  \odot : {\rm M} (n \times m, A) \otimes {\rm M} (k \times l, A) \rightarrow
  {\rm M} (nk \times ml, A) 
\end{equation}
 defined by
\begin{equation}
   ( M \odot N)_{\alpha \iota, \beta \kappa} = M_{\alpha, \beta} \, 
  N_{\iota, \kappa}.  
\end{equation} 

 Now let $V$ be a $n$-dimensional complex vector space and let 
 $R : \Omega \rightarrow {\rm End \,} (V) \otimes {\rm End \,}(V)$ 
 be a holomorphic map defined on an open domain $\Omega \subset \C$ such that 
 $ 1 \in \Omega$ and $R(1) = 1$.
 Choosing a basis $(x_1,...,x_n)$ of $V$ we can regard $R$ as a $n \times n$
 matrix with entries in ${\cal O} (\Omega)$. The basis $(x_1,...,x_n)$
 naturally induces a basis $(t_{\iota}^{\kappa})$ of the coalgebra 
 $C = ( \End \, V )'$. Furthermore, it gives rise to the matrix 
 $T \in {\rm M} (n \times n , {\rm T}C ) $ the entries
 of which are the $t_{\iota}^{\kappa}$.  
 The tensor algebra ${\rm T} C$ shall carry the inductive limit 
 topology of all finite dimensional subspaces. 
 Thus it is a nuclear bialgebra with coproduct 
\begin{equation}
  \Delta (t_{\iota}^{\kappa}) = \sum \limits_{1 \leq \alpha \leq n} \,
  t_{\iota}^{\alpha} \otimes t_{\alpha}^{\kappa}.
\end{equation}
and counit
\begin{equation}
  \epsilon (t_{\iota}^{\kappa}) = \delta_{\iota}^{\kappa}.
\end{equation}
 We now construct the 
 {\bf quantum matrix algebra} or {\bf FRT-bialgebra} $A(R)$ as the quotient of
 ${\rm T} C$ modulo the closed ideal $I$ generated by the 
 coefficients of the matrix
\begin{equation}
 ( T \odot 1 )\, (1 \odot T)\,  R \:  - \: R \, (1 \odot T) \, (T \odot 1).
\end{equation}
\end{rPar}
 By {\sc Larson and Tauber} \cite{LarTow:QGQLA}  or 
 {\sc Manin} \cite{Man:QGNG}  $I$ is a biideal, so
 $A(R)$ becomes a nuclear bialgebra with a continuous coaction on 
 ${\cal O} (\Omega) \otimes V$ given by
\begin{equation}
  \Psi_{V,R} : {\cal O} (\Omega) \otimes V \rightarrow A(R) \otimes V, \quad
  x_{\iota} \mapsto t_{\iota}^{\kappa} \otimes x_{\kappa}.
\end{equation}
 Using Theorem and Definition 3.1 of 
 {\sc Larson and Tauber} \cite{LarTow:QGQLA} it follows that $A(R)$ is
 the initial object with respect to the following properties: \\[3mm]
 (FRT1) ${\cal O} (\Omega) \otimes V$ is a topological left comodule over
 $A(R)$ via the structure map $\Psi_{V,R}$. \\[2mm]
 (FRT2) $R$ is an $A(R)$-comodule map with respect to the natural structure of 
 a topological\\ \mbox{ }\hspace{12mm} left $A(R)$-comodule on 
 $  {\cal O} (\Omega) \otimes V \otimes V$. \vspace{3mm}
 
 If $R$ is a Yang-Baxter operator, i.e.~$R$ is 
 non-degenerate and fulfills the quantum Yang-Baxter equation (QYBE)
\begin{equation}
 R_{12} \circ R_{13} \circ R_{23} \: = \: R_{23} \circ R_{13} \circ R_{12},
\end{equation}
 then Theorem and Definition 5.1 of {\sc Larson and Tauber} \cite{LarTow:QGQLA}
 imply that $A(R)$ is topologically coquasitriangular.
 The braiding form  $< \: | \: > : \, A(R) \otimes A(R) \rightarrow \C$ 
 can furthermore be calculated by the relations 
\begin{equation}
 < t_{\iota}^i | t_{\kappa}^j > \: = \: R_{\iota \kappa}^{ji}.
\end{equation}

\section{Holomorphic deformation of locally convex algebras}
\label{HoDeLoCoAl}
In our approach to deformation theory we replace the ring $\C [[ h ]]$
used in the formal deformation theory of {\sc Gerstenhaber} 
\cite{Ger:DRA,GerSch:ACDT,BonFlaGerPin:HGSQGSDRPD} by the ring ${\cal O} (\Omega)$ of holomorphic 
functions on an open domain $\Omega \subset \C^n$. 
Note that in the case of $\Omega = \C$
we have the following relations and inclusions:
\begin{eqnarray}
\begin{array}{ccccc}
  \C [h] & \longmapsto & \C [[ h ]] & \longmapsto & \C [h] \\
  \cap & & \cup & & \cap\\
  {\cal O} (\C) & \longmapsto & {\cal O} \{ z \} & \longmapsto & {\cal O}(\C) \\
  \cap & & \cup & & \cap \\
  \C [[ h ]] & \longmapsto &  \C [h] & \longmapsto & \C [[ h ]]   ,
\end{array}
\end{eqnarray}
where the arrows denote passing to the topological dual, and where 
${\cal O} \{ z \}$ denotes the algebra of convergent power series in
$z$ with the natural locally convex inductive limit topology. 

In the following let $\Omega$ be an open domain in $\C^n$ and let ${\cal O} =
{\cal O} (\Omega)$ be the nuclear Fr\'echet algebra of holomorphic functions 
on $\Omega$. For a complete locally convex Hausdorff space $E$ let  
\begin{equation}
  E_{\Omega} \: = \: {\cal O} (\Omega , E) 
\end{equation}
be the space of holomorphic $E$-valued functions $f : \Omega \rightarrow E$
equipped with the compact open topology. Every $E$ can locally be represented
by a convergent power series 
$f(z) = \sum_{\alpha \in \N^n} \, f_{\alpha} z^{\alpha}$, 
$z \in \Omega$ where $f_{\alpha} \in E$ (cf.~Appendix). $E_{\Omega}$ 
is a complete locally convex Hausdorff space which is isomorphic to the 
completion of the tensor product ${\cal O} (\Omega) \otimes E$ with respect to the $\pi$-topology 
(cf.~Appendix):
\begin{equation}
\label{EqFuTePr}
  E_{\Omega} \: = \: {\cal O} (\Omega , E) \: \cong \: {\cal O} (\Omega) 
  \hat{\otimes} E.
\end{equation}
Hence, $E_{\Omega}$ is nuclear (resp.~Fr\'echet), if $E$ is nuclear 
(resp.~Fr\'echet).
The pointwise multiplication 
${\cal O} (\Omega) \times E_{\Omega} \rightarrow E_{\Omega},$ 
$ ( \lambda, f) \mapsto \lambda f $ 
is continuous and thus defines on $E_{\Omega}$ the structure of a 
{\bf topological ${\cal O} (\Omega)$-module}. 
Because of Eq.~(\ref{EqFuTePr}) we sometimes call $E_{\Omega}$ a 
{\bf topologically free ${\cal O} (\Omega)\,$}-{\bf module}.

For a locally convex algebra $A$ the space 
\begin{equation}
  A_{\Omega} \: = \: {\cal O} (\Omega , A) \: \cong {\cal O} (\Omega  ) 
  \hat{\otimes} A 
\end{equation}
is a locally convex algebra and a topological module over ${\cal O} (\Omega)$.
The algebra  structure on $A_{\Omega}$ which we call the {\bf constant algebra
structure} is given by pointwise multiplication and trivial extension of the
unit: If $\mu : A \hattimes A \rightarrow A$ and $\eta : \C \rightarrow A$ 
are multiplication and unit of $A$ respecticely, the corresponding 
structure maps of the constant algebra structure on $A_{\Omega}$ are 
given by
\begin{eqnarray}
  \mu = \mu_{\Omega} : & A_{\Omega} \times A_{\Omega} \rightarrow A_{\Omega}, &
  (f,g) \mapsto
  \big( \Omega \ni z \mapsto \mu_{\Omega} (f,g) (z) = \mu( f(z) , g(z)) \in A
  \big)  
\end{eqnarray}
  and 
\begin{eqnarray}
  \eta = \eta_{\Omega} : & \C \rightarrow A_{\Omega}, & \lambda \mapsto 
  \lambda \otimes \eta (1).
\end{eqnarray}
Evidently the mappings
$\mu_{\Omega}$ and $\eta_{\Omega}$ are continuous and fulfill the axioms of
associativity and unit. They are in fact holomorphic in an obvious sense. 
Under the isomorphism 
${\cal O} (\Omega , A) \: \cong {\cal O} (\Omega  ) \hattimes A$
the product $\mu_{\Omega}$ is given by
\begin{equation}
  A_{\Omega} \otimes  A_{\Omega} \ni (f \otimes a) \otimes (g \otimes b) 
  \mapsto (fg \otimes ab) \in  A_{\Omega}.
\end{equation}

In case $A$ is a complete locally m-convex algebra we have an obvious 
functional calculus on $A_{\Omega}$, i.e.~for any holomorphic function 
$f = \sum_{\alpha \in \N^n} \, f_{\alpha} z^{\alpha}$ in ${\cal O}(\C)$  
and any element $a \in A$ there is a unique 
$f(a)= \sum_{\alpha \in \N^n} \, f_{\alpha} a^{\alpha} \in A$
(cf.~Appendix).
\vspace{0.2cm}

In deformation theory we are now interested in algebra structures
on $A_{\Omega}$ different from the constant one. Let us describe this in
more detail.
Define a {\bf topological ${\cal O} (\Omega)\,$-algebra structure}
on $E_{\Omega}$ to be an algebra structure on 
$E_{\Omega} = {\cal O} (\Omega,E)$ 
given by a continuous ${\cal O} (\Omega)\,$-bilinear map
\begin{equation}
  \tilde{\mu} : E_{\Omega} \times E_{\Omega} \longrightarrow E_{\Omega}
\end{equation}
 fulfilling the associativity constraint and a continuous 
${\cal O}(\Omega)$-linear unit
\begin{equation}
  \tilde{\eta} : {\cal O}(\Omega) \longrightarrow E_{\Omega}.
\end{equation}
  We often denote this algebra 
  $(E_{\Omega}, \tilde{\mu}, \tilde{\eta})$ by $\tilde{E}$. 
  $\tilde{\eta}$ is determined by $\tilde{\eta}(1) \in E_{\Omega}$.
  
  Two such algebra structures $(\tilde{\mu}, \tilde{\eta})$ and 
  $(\check{\mu}, \check{\eta})$ on $E_{\Omega}$ are called {\bf equivalent}
  if  there exists a ${\cal O} (\Omega)$-linear isomorphism
  $\varphi : E_{\Omega} \longrightarrow E_{\Omega}$  (of locally convex spaces)
  such that the relations
\begin{eqnarray}
  \varphi \circ \check{\mu} & = & \tilde{\mu} \circ (\varphi \times 
  \varphi), \\
  \varphi \circ \check{\eta} & = & \tilde{\eta}
\end{eqnarray}
  hold.

 The multiplication $\tilde{\mu}$ of an ${\cal O} (\Omega)$-algebra 
 $(E_{\Omega}, \tilde{\mu}, \tilde{\eta})$ can be regarded to be a
 continuous ${\cal O} (\Omega)$-linear map $\tilde{\mu} : E_{\Omega}
 \hattimes_{{\cal O} (\Omega)} E_{\Omega} \rightarrow  E_{\Omega}$, where
 $E_{{\Omega} } \hattimes_{{\cal O} (\Omega)} E_{\Omega}$ is isomorphic (as a
 topological ${\cal O} (\Omega)$-module) to ${\cal O} (\Omega , E {\hattimes}
 E)$ (cf.~Appendix). 
 Consequently $\tilde{\mu}$ can be described as a
 holomorphic map $\tilde{\mu} : \Omega \rightarrow {\cal L} ( 
 E {\hattimes} E, E)$, where ${\cal L} (E,F)$ is the space of continuous 
 linear maps from $E$ to $F$ equipped with the topology of uniform
 convergence on bounded sets of $E$.
 Thus around any $z_0 \in \Omega$ the product map $\tilde{\mu}$ can locally
 be expanded in the form
\begin{equation}
\label{ExMu}
  \tilde{\mu} (z) \: = \: \sum \limits_{\alpha \in \N^n} \, \mu_{\alpha} 
  (z - z_0)^{\alpha}, 
\end{equation}
where $z \in \Omega$ is close enough to $z_0$ and the $\mu_{\alpha}$
are continuous bilinear mappings on $E$ with values  in
$E$ and which depend on the base point $z_0$. 
Similarly $\tilde{\eta} : \C \rightarrow E_{\Omega}$ with $\tilde{\eta}
(\lambda) = \lambda \, \tilde{\eta} (1)$ can be viewed as the map
$\tilde{\eta} (1) : \Omega \rightarrow E$.
In case $E = A$ is a locally convex algebra the constant algebraic structure 
on $A_{\Omega}$ is then given by $\tilde{\mu} (z) = \mu$ and 
$\tilde{\eta} (\lambda) (z) = \eta (\lambda)$, $\lambda \in \C$, 
$z \in \Omega$.
 
Before giving the definition of a deformation let us mention that
in the following $\m_z$ denotes the maximal ideal of ${\cal O} 
(\Omega)$ at a point $z \in \Omega$, i.e.~$\m_z$ is the ideal
of all holomorphic functions on $\Omega$ vanishing at $z$.
\begin{Definition}
  A {\bf topologically free holomorphic deformation} 
  of a locally convex algebra $(A,\mu,\eta)$
  over a complex domain $\Omega$ is an ${\cal O} (\Omega)\,$-algebra 
  structure $(\tilde{\mu},\tilde{\eta})$ on $A_{\Omega}$ such that for a
  distinguished point $* \in \Omega$
  the quotient ${\cal O} (\Omega)$-module $A_{\Omega} / \m_* A_{\Omega}$ 
  is isomorphic to $A$ as a locally convex algebra.
  Equivalently $\mu_0 = \mu$ in the expansion (\ref{ExMu}) and 
  $\tilde{\eta}(1)(*) = \eta(1)$. 
  
  The deformation is called {\bf trivial} if 
  $(A_{\Omega}, \tilde{\mu}, \tilde{\eta})$ is equivalent to the constant 
  algebra structure on $A_{\Omega}$. The distinguished point $*$ is called
  the {\bf base point} of the deformation.
\end{Definition}

In the sequel we shall call a topologically free  holomorphic deformation 
simply a holomorphic deformation, although there exist more general 
schemes of holomorphic deformations 
(see {\sc Pflaum} \cite{Pfl:LADQ,Pfl:NCDQNOQCB}).
However, these are not needed for our purposes of deformation of
Hopf algebra structures.

An important advantage of holomorphic deformations 
in comparison to formal deformations lies in the fact that
for every parameter $ z \in \Omega$ one receives a {\bf concrete}
deformed algebra structure on the underlying linear space of the original 
algebra $A$: Simply take -- for any value $z \in \Omega$ -- 
$\tilde{\mu} (z) \in {\cal  L} ( A \hattimes A , A)$
as the deformed multiplication and $\tilde{\eta} (z) \in {\cal L} (\C, A)$ 
as the deformed unit.
Then for $a,b \in A$
\begin{equation}
  a \, *_z \, b \: = \: \tilde{\mu} (z) \, (a,b) \in A 
\end{equation}
is the new product of $a$ and $b$ and 
\begin{equation}
   e_z \: = \: \tilde{\eta} (z) \, (1) \in A
\end{equation}
the new unit. Both, the new product and the new unit are contained in the 
original space $A$ and not only in a nontrivial extension of $A$.

Of course, the question now is, whether there exist interesting and
nontrivial holomorphic deformations of locally convex algebras. 
The main purpose of this paper is to answer this question in the positive.
\begin{rExample}
\begin{enumerate}
\item
  {\bf Quantum vector spaces} (cf.~{\sc Faddeev, Reshetikhin, Takhtajan}, 
  \cite{FadResTak:QLGLA} {\sc Manin} \cite{Man:QGNG}). 
  Let $\Omega = \C^* $ be the set of all
  nonzero complex numbers and consider the $n$-dimensional complex 
  vector space $V = \C^n$ with the canonical basis $(x_1,...,x_n)$.
  Then  construct the tensor algebra ${\rm T}V$ of $V$ or in other words the 
  free $\C$-algebra in $n$ generators. By completion with respect to
  the projective limit topology of all Fr\'echet representations ${\rm T}V$ 
  becomes a nuclear locally m-convex Fr\'echet algebra $\check{\rm T}V$. 
  The functions
\begin{equation}
\label{CoQuVeSp}
  f_{\iota \kappa} : \: \Omega \rightarrow {\rm T} V , \quad z \mapsto 
  x_{\iota} \, x_{\kappa} - z \, x_{\kappa} \,  x_{\iota} \quad
  1 \leq \iota < \kappa \leq n
\end{equation}
  then generate a unique closed ideal $I$ in the nuclear algebra
  ${\cal O} (\Omega) \hattimes \check{\rm T} V$. The quotient
  ${\cal O} (\C^n_q ) = {\cal O} (\Omega) \hattimes \check{\rm T} V/ I$ 
  is called the 
  {\bf algebra of entire functions on the quantum $n$-vector space}.
  It comprises a holomorphic deformation of the algebra ${\cal O} (\C^n)$
  of entire functions on $\C^n$: 
  The homomorphism ${\cal O} (\C^n_q) \rightarrow {\cal O} (\C^n)$
  defined by $ [f \otimes  x_{\iota} ] \mapsto f(1) x_{\iota} $ is
  well-defined, surjective and has kernel $\m_1 {\cal O} (\C^n_q)$.
  Therefore 
  ${\cal O} (\C^n_q) / \m_1 {\cal O} (\C^n_q) \cong {\cal O} (\C^n)$    
  holds. Because of the relations (\ref{CoQuVeSp}) the 
  ${\cal O} (\C^*)$-linear combinations of the  family 
  $(x_1^{m_1},...,x_n^{m_n})_{m_1,...,m_n \in \N}$ are dense in 
  ${\cal O} (\C_q^n)$. Since 
  $(x_1^{m_1},...,x_n^{m_n})_{m_1,...,m_n \in \N}$  furthermore
  is free over ${\cal O} (\C^*)$, ${\cal O} (\C_q^n)$
  is isomorphic to ${\cal O} (\C^*)\hattimes {\cal O} (\C^n) \cong {\cal O}
  (\C^* \times \C^n_q)$ as a nuclear space. This proves the claim.  
  
  Alternatively one could give ${\rm T} V$ the inductive topology
  of all finite dimensional subspaces. Then ${\rm T} V$ is already a strictly
  nuclear algebra. By the same procedure as above but now applied to
  ${\rm T} V$ one is lead to the algebra ${\cal P} (\C_q^n)$ of 
  {\bf polynomial functions on the quantum $n$-vector space}. 
  ${\cal P} (\C_q^n)$ comprises a deformation of the algebra $\C [x_1,...,x_n]$
  of polynomials in $n$ complex variables.
\item
  {\bf Quantum exterior algebra} (cf.~{\sc Wess, Zumino} \cite{WesZum:CDCQH},
  {\sc Manin} \cite{Man:QGNG}).
  In the spirit of the preceeding example
  it is also possible to deform the exterior algebra on $\C^n$. 
  Let $V'$ be the dual of $V$, $(\xi_1,...,\xi_n)$ 
  the dual basis of $(x_1,...,x_n)$ and let ${\rm T} V'$ be given the finest 
  locally convex topology.  Then consider the closed ideal 
  $J \subset {\cal O} (\C^*) \hattimes {\rm T} V' $
  generated by the relations
\begin{equation}
\label{CoQuDiFo}
  \xi_{\iota}^2 = 0,  \quad
  \xi_{\iota} \, \xi_{\kappa} \: = \: - z^{-1} \, \xi_{\kappa} \,  
  \xi_{\iota} \quad  1 \leq \iota < \kappa \leq n.
\end{equation} 
  The corresponding quotient $\Lambda (\C_q^n) = {\rm T} V' / J $ 
  is the {\bf exterior algebra of the quantum $n$-vector space}.
  Exactly like above it is shown that $\Lambda (\C^n_q )$ is a holomorphic
  deformation over $\C^*$ of the exterior algebra $\Lambda (\C^n)$.
  Note that unlike ${\cal O} (\C_q^n)$ and ${\cal P}  (\C_q^n)$
  the algebra $\Lambda (\C^n)$ is finite dimensional. 
  The tensor product algebra ${\cal O} (\C^n_q) \otimes
  \Lambda (\C^n)$ can be interpreted as the algebra of 
  {\bf entire  holomorphic quantum differential forms on $\C^n_q$},
  the tensor product ${\cal P} (\C_q^n)\otimes \Lambda (\C^n) $
  as the algebra of 
  {\bf algebraic quantum differential forms on $\C^n_q$}.
\end{enumerate}
\end{rExample}
\begin{rPar} {\bf Holomorphic Quantization.}
 In mathematical physics one is not only interested in deformations of
 algebras but also in their quantization. More precisely we want to
 quantize so-called {\bf locally convex Poisson algebras}, i.e.~locally
 convex $\K$-algebras $A$ on which a continuous bilinear bracket 
 $\{ \hspace{0.5em} , \hspace{0.5em} \} : A \times A \rightarrow A$
 is defined which is antisymmetric, fulfills the Jacobi identity
 and the Leibniz rule. 
 A {\bf topologically free holomorphic quantization} of $A$ over the domain 
 $\Omega \subset \C$ with base point $0 \in \Omega$ then is a 
 topologically free holomorphic deformation 
 $(A_{\Omega}, \tilde{\mu}, \tilde{\eta})$ of $A$ such that the relation
\begin{equation}
  \tilde{\mu} (f, g) - \tilde{\mu} (g,f) \: = \: - i \, z \, \{ f , g \} 
  \: + \: o (z^2)
\end{equation}  
  holds for all $f,g \in  A $ and $z \in \Omega$.
\begin{rRemark}
\label{HQFFQ}
  Suppose  $A$ to be a locally convex Poisson algebra having a holomorphic 
  deformation $(\tilde{A}, \tilde{\mu})$. 
  If now $\tilde{A} \hattimes \C[[\hbar]]$ is a formal quantization of $A$, 
  then $\tilde{A} $ obviously is a holomorphic quantization of $A$.
\end{rRemark}
\end{rPar}

Let us give an important example of a holomorphic quantization. Consider the
space $\symi (\R^n \times \R^n)$ of classical symbols on $\R^n$
together with the topology of asymptotic convergence (cf.~{\sc H\"ormander}
\cite{Hor:ALPDO} or {\sc Pflaum} \cite{Pfl:NCDQNOQCB} for details on symbol
spaces). The Poisson bracket on $\R^{2n}$ naturally induces a Poisson 
bracket on the symbol space $\symi (\R^n \times \R^n)$. We want to quantize
this Poisson space. By the asymptotic expansion
\begin{equation}
\label{MoPr}
  ( a *_{\hbar} b )
  \: (x , \xi ) \sim \sum \limits_{\alpha, \beta  \in \N^n} \, 
  (-1)^{|\beta|} \, \left( \frac{- i \hbar}{2} \right)^{|\alpha +\beta| }\, 
  \left( \dera{\xi} \derb{x} a (x, \xi) \right) \,
  \left( \derb{\xi} \dera{x} b (x, \xi) \right)
\end{equation} 
we define the {\bf Moyal-product} $a *_{\hbar} b $ of two symbols
$a,b \in \symmi (\R^n \times \R^n)$.
Unfortunately -- as asymptotic expansions are only unique up to smoothing 
symbols $\in \symmi (\R^n \times \R^n)$ -- Eq.~(\ref{MoPr}) does not define
an associative product on $\symi (\R^n \times \R^n)$. 
But the Moyal product gives a holomorphic deformation of the quotient algebra
$\symi / \symmi (\R^n \times \R^n)$ carrying the topology of asymptotic
convergence. As is well-known the Moyal product induces a formal
quantization of the Poisson space ${\cal C}^{\infty} (\R^n \times \R^n)$, 
hence, of $\symi (\R^n \times \R^n)$. As the Poisson bracket on 
$\symi (\R^n \times \R^n)$ can be pushed down to the quotient 
$ \symi / \symmi ( \R^n \times \R^n )$, the above proposition entails
that the Moyal-product induces a holomorphic quantization of
$ \symi / \symmi ( \R^n \times \R^n )$. 

Note that in the same way the Wick quantization procedure leads to a 
holomophic quantization. See also {\sc Pflaum} \cite{Pfl:NCDQNOQCB} for more
general quantization schemes which fit into the above holomorphic picture.

\section{Holomorphic deformation of nuclear Hopf algebras}
\label{HoDeNuHoAl}
The concept of a (topologically free) holomorphic deformation can easily 
be transferred to the case of deformations of nuclear coalgebra structures, 
bialgebra and Hopf algebra structures, as well as Lie bialgebra structures. 
In particular a {\bf (topologically free) holomorphic deformation}
of a nuclear Hopf algebra $H$ with multiplication $\mu$, comultiplication 
$\Delta$, unit $\eta$, counit $\varepsilon$, and antipode $S$ is given by 
the following data:
\begin{eqnarray}
\begin{array}{lclclclcl}
	\tilde{\mu} & \!\!\! \in \!\!\! & {\cal O} 
	(\Omega,{\cal L}( H \hat{\otimes} H,H)
	& \cong & {\cal L}_{{\cal O} (\Omega)} (H_{\Omega} \hat
	{\otimes}_{{\cal O} (\Omega)} H_{\Omega}, H_{\Omega})	
	& \mbox{with} &  
	\tilde{\mu} (*) &\! =\!& {\mu},\\

	\tilde{\Delta} & \!\!\! \in \!\!\!  & {\cal O} (\Omega,{\cal L}(H, 
        H \hat{\otimes} H)
	& \cong & {\cal L}_{{\cal O} (\Omega)} (H_{\Omega} , H_{\Omega} \hat
	{\otimes}_{{\cal O} (\Omega)} H_{\Omega})	
	& \mbox{with} & 
	\tilde{\Delta}(*) \!& = \!& {\Delta},\\

	\tilde{\eta} & \!\!\! \in \!\!\! & {\cal O} (\Omega,{\cal L}(\C, H))
	& \cong & {\cal L}_{{\cal O} (\Omega)} ({\cal O} (\Omega), H_{\Omega}) 
	& \mbox{with} &
	\tilde{\eta} (*) \!& = \!& {\eta},\\

	\tilde{\epsilon} & \!\!\! \in \!\!\! & {\cal O} (\Omega,{\cal L}( H,\C))
	& \cong &{\cal L}_{{\cal O} (\Omega)}  (H_{\Omega},{\cal O} (\Omega))
        & \mbox{with} & 
	\tilde{\epsilon} (*) \!& = \!& {\epsilon},\\

	\tilde{S} & \!\!\! \in \!\!\! & {\cal O} (\Omega, {\cal L}(H, H))
	& \cong & {\cal L}_{{\cal O} (\Omega)} (H_{\Omega}, H_{\Omega})
	& \mbox{with} & 
	\tilde{S} (*) \!& = \!& S, \\
\end{array}
\end{eqnarray}
such that $(H_{\Omega},\tilde{\mu},\tilde{\Delta}, \tilde{\eta}, 
\tilde{\epsilon},\tilde{S})$ is a nuclear Hopf algebra.
(The isomorphisms in the above formulas are explained in the Appendix.) 
Two such Hopf algebra deformations are {\bf equivalent} if the corresponding
Hopf algebra structures on $H_{\Omega}$ are isomorphic.
Any holomorphic deformation of a nuclear Hopf algebra turns out to be 
equivalent to a holomorphic deformation with constant unit and counit, 
i.e. with $\tilde{\eta} (z) = {\eta}$ and  $\tilde{\epsilon} (z) = {\epsilon}$
for all $z \in \Omega$.

Our first result concerns the dual of the holomorphic deformation $\tilde{H}$
of a strictly nuclear Hopf algebra $H$: 
Transposing all the structure maps one gets a 
holomorphic deformation of the dual nuclear Hopf algebra $H'$. 
The underlying topological ${\cal O} (\Omega)$-module is 
\begin{equation}
 {H_{\Omega}}' \: = \: {\cal L}_{{\cal O} (\Omega)} 
 (H_{\Omega}, {\cal O} (\Omega)) =
 \{\phi \,:\: H_{\Omega} \to {\cal O} (\Omega) \;|\; \phi \;
 \mbox{is continuous and} \;{\cal O} (\Omega)\mbox{-linear} \}
\end{equation}
with the topology of uniform convergence on bounded (or equivalently 
compact) sets of $H_{\Omega} \cong {\cal O} (\Omega) \hat{\otimes} H$.
This space is isomorphic as an ${\cal O} (\Omega)$-module to
\begin{equation}
  {\cal O} (\Omega ) \hat{\otimes} H' \: \cong
  {\cal O} (\Omega , H' ) \: = \: H'_{\Omega} .
\end{equation}
The structure maps on $H_{\Omega}'$ are
\begin{equation}
  \tilde{\Delta}' \: = \: 
  ( \tilde{\mu} )^* : \: H_{\Omega}' \longrightarrow  
  ( H_{\Omega} \hat{\otimes} H_{\Omega})'  \cong 
  H_{\Omega}' \hattimes H_{\Omega}' 
  \cong {\cal O} (\Omega, H' \hat{\otimes} H'),
\end{equation}
hence 
$\tilde{\Delta}' 
\in {\cal O} (\Omega , {\cal L} (H' , H' \hat{\otimes} H'))$,
\begin{equation}
  \tilde{\mu}' \: = \: ( \tilde{\Delta} )^* : \:
  ( H_{\Omega}  \hat{\otimes} H_{\Omega})' \cong
  H_{\Omega}' \otimes H_{\Omega}' \rightarrow H'_{\Omega} ,
\end{equation}
hence 
$\tilde{\mu}' \in {\cal O} ( \Omega , {\cal L} (H' \hat{\otimes} H' , H' ))$,
and similarly 
\begin{eqnarray}
  \tilde{S}' \: = \: \tilde{S}^*, & 
  \tilde{\eta}' \: = \: \tilde{\varepsilon}^*, &   
  \tilde{\varepsilon}' \: = \: \tilde{\eta}^*.
\end{eqnarray}
The following duality theorem is now obvious.
\begin{Theorem}
  Let $\tilde{H} = ( H_{\Omega}, \tilde{\mu},  \tilde{\Delta}, 
  \tilde{\varepsilon},\tilde{\eta}, \tilde{S} )$ 
  be a holomorphic deformation of the strictly nuclear 
  Hopf algebra $(H, \mu , \Delta, \varepsilon, \eta , S)$.
  Then the ${\cal O} (\Omega)$-dual $H_{\Omega}'$ of $H_{\Omega}$ with the 
  structure maps $\tilde{\mu}'$, $\tilde{\Delta}'$, 
  $\tilde{\eta}'$, $\tilde{\varepsilon}'$ and $\tilde{S}'$
  is  a holomorphic deformation of the nuclear Hopf algebra
  $H'$. Moreover, $\tilde{H}$ and $\check{H}$ are equivalent if and only
  if $\tilde{H}'$ and $\check{H}'$ are equivalent.
  In addition, the bidual $(\tilde{H}')'$ is canonically
  isomorphic to $\tilde{H}$.
\end{Theorem}
The above theorem is a direct generalization of the corresponding result
for the formal case (cf.~\cite{BonFlaGerPin:HGSQGSDRPD}). Similarly, the next
fact about the construction of a twisting matrix
carries over directly to the holomorphic case 
(cf.~\cite{BonFlaGerPin:HGSQGSDRPD}, Proposition 4.2.4).
\begin{Theorem}
\label{TwMaUnEnAl}
  Let $G$ be a compact group and let $H$ denote the strictly nuclear Hopf 
  algebra ${\cal R} (G)' \cong \tilde {\cal U} \g^{\C}$ respectively
  ${\cal D} (G) = {\cal E} (G)'$ (cf.~\ref{UnEnAl} resp.~\ref{FuDiSpLiGr}). 
  Let $(H_{\Omega} , \tilde{\Delta})$ be a holomorphic deformation of the 
  nuclear bialgebra $H$ leaving invariant the algebra product on $H$. Then
  there exists ${\cal F} \in \left( H \hattimes H \right)_{\Omega}$ such 
  that $\tilde{\Delta} \, {\cal F} = {\cal F} \, \Delta$.  
\end{Theorem}
\begin{Proof}
Because of $\tilde{\Delta} \in {\cal L}_{ {\cal O} ( \Omega ) } ( H_{\Omega }, 
H_{\Omega } \otimes_{ {\cal O} (\Omega ) } H_{\Omega } ) \cong 
{\cal O} ( \Omega , {\cal L} ( H, H \hattimes H))$ we get a ``twisting
matrix'' ${\cal F}  \in {\cal O} (\Omega , H \hattimes H ) \cong 
(H \hattimes H)_{{\cal O} (\Omega)}$ by the following integral:
\begin{equation}
  {\cal F} \: = \: \int \limits_G \, \tilde{\Delta} (g ) \, \left( \Delta (g)
  \right)^{-1} \, d \mu (g),
\end{equation}
  where $\mu$ is the Haar measure on $G$ and $g \in G$ stands for $\delta_g$, 
  the Dirac distribution with support in $g$.
  By left- and right-invariance of the Haar measure
\begin{eqnarray}
  {\cal F} \, \Delta (h) & = & \int \limits_G \, \tilde{\Delta} (hg) \,
  ( \Delta (hg))^{-1} \, \Delta (h) \, d \mu (hg) \: =  \\ \nonumber
  & = & \int \limits_G \tilde{\Delta} (h) \, \tilde{ \Delta} (g) 
  \, (\Delta (g) )^{-1} \, d\mu (g)  \: = \: \tilde{\Delta} (h) \, {\cal F}
\end{eqnarray}
  holds. Hence the claim follows from the fact that the Dirac distributions
  $g = \delta_g$ lie densely in ${\cal R} (G)' (cf.~\ref{MaCoGrRel})$. 
\end{Proof}
\begin{rRemark}
 We expect in all important examples ${\cal F}(z)$ to be invertible at least 
 in a neighborhood of the base point. In that case we get a dual version of
 Theorem \ref{TwMaUnEnAl} which leads to a holomorphic coquasitriangular
 Hopf algebra deformation of ${\cal R} (G)$.
\end{rRemark}
\begin{Theorem}
  Let $(H,\mu,\eta,\Delta,\varepsilon,S)$ be a cocommutative 
  nuclear Hopf algebra and
  let ${\cal F}: \Omega \rightarrow H \hattimes H$ be 
  a {\bf twisting matrix} for $H$, i.e.~a continuous map
  on an open domain $\Omega \subset \C$ such that  
  the following conditions hold:   
\begin{enumerate}
\item
  ${\cal F}(z)$ is an invertible element of $H \hattimes H$ for every
  $z \in \Omega$.
\item
  ${\cal F}$ fulfills the equations
\begin{eqnarray}
  {\cal F}_{12} \: ( \Delta \otimes 1 ) \: {\cal F} & = & 
  {\cal F}_{23} \: ( 1 \otimes \Delta ) \: {\cal F} \\
  (\varepsilon \otimes 1) {\cal F} \: = & 1 & = \:
  ( 1 \otimes \varepsilon ) {\cal F},
\end{eqnarray}
  where ${\cal F}_{12} $ and ${\cal F}_{23} $ 
  are defined as usual.
\end{enumerate}
 Then ${\cal F}$ induces a new nuclear Hopf algebra structure 
 $(H_{\Omega},\mu^{\cal F},\eta^{\cal F},\Delta^{\cal F},\varepsilon^{\cal F},
 S^{\cal F})$ on $H_{\Omega} = {\cal O} (\Omega) \hattimes H$ by defining 
 $\mu^{\cal F} = \mu $, $\eta^{\cal F} =  \eta $,
 $\Delta^{\cal F} (z \otimes h) = {\cal F} (z) \, \Delta (h) {\cal F}^{-1} (z)$, 
 $\varepsilon^{\cal F} = \varepsilon$ and $S^{\cal F} (z \otimes h) =
 v (z) \, S (h) \, v^{-1} (z)$. Hereby is $z \in {\cal O} (\Omega)$,
 $h \in H$ and $v = \mu ( 1 \otimes S ) {\cal F}$.
 $(H_{\Omega}, \mu, \eta, \Delta^{\cal F} , \varepsilon, S^{\cal F})$
 is topologically triangular with universal R-matrix 
 ${\cal R} = {\cal F}_{21} \: {\cal F}^{-1}$.

 Now assume $H$ to be topologically quasitriangular with universal R-matrix 
 ${\cal R}$ and, additionally to the above, the following:
\begin{enumerate}
\setcounter{enumi}{2}
\item 
  The quantum Yang-Baxter equation (QYBE) is fulfilled:
\begin{equation}
  {\cal F}_{12} \circ {\cal F}_{13} \circ {\cal F}_{23} \: = \: 
  {\cal F}_{23} \circ {\cal F}_{13} \circ {\cal F}_{12}.
\end{equation} 
\item 
  ${\cal F}_{21} = {\cal F}^{-1}$.
\end{enumerate} 
  Then $(H_{\Omega},\mu, \eta, \Delta^{\cal F},\varepsilon, S^{\cal F})$ 
  is a   topologically quasitriangular nuclear Hopf algebra with
  universal R-matrix 
  ${\cal R}^{\cal F} = {\cal F}^{-1} {\cal R} {\cal F}^{-1}$. 

  If the twisting matrix ${\cal F}$ fulfills ${\cal F} (*) = 1 \otimes 1$ 
  for a point $* \in \Omega$,
  then in both of the above cases 
  $(H_{\Omega},\mu,\eta, \Delta^{\cal F},\varepsilon,S^{\cal F})$
  comprises a topologically free holomorphic deformation over $\Omega$
  of the Hopf algebra $H$ with base point $*$. 
\end{Theorem}
\begin{Proof}
The proof of the theorem can be taken almost literally from the corresponding 
one in the formal case. Confer for example {\sc Chari, Presley} 
\cite{ChaPre:GQG}.
\end{Proof}

The purpose of the following considerations is to show how one can construct 
under certain conditions a holomorphic deformation of an algebra 
(resp.~bialgebra or Hopf algebra) out of a formal one. We will carry out
the details only for the algebra case; the ones for bialgebras and 
Hopf algebras are analogous.

So let $A$ be a finitely generated $\C$-algebra and assume that there exists
a formal deformation $A_{\hbar} \cong (A[[\hbar]], \mu,\eta)$ of $A$.
In other words $\mu$ is a $\C[[\hbar]]$-bilinear multiplication map on 
$A[[\hbar]]$ and $\eta$ a unit such that 
$A[[\hbar]]/ \hbar A[[\hbar]] \cong A$.
We now choose a finite dimensional vector space $V$ and a surjective 
homomorphism ${\rm T}V = \bigoplus_{n \in \N} V^{\otimes n} \rightarrow A$.
According to \ref{PrLiToHoAl} this gives rise to a topological presentation
$\hat{\rm T} V \rightarrow \hat{A} $, where $\hat{A}$ 
and $\hat{\rm T}V$ are the completions of $A$ resp.~${\rm T}V$ with respect to 
the topology of all finite dimensional representations
(resp.~all nuclear Fr\'echet representations).
Using the universal property of the tensor algebra ${\rm T}V$ there exists
a unique morphism of $\C[[\hbar]]$-algebras 
$\pi: \, {\rm T} V [[\hbar]] \rightarrow A_{\hbar}$  such that 
${\rm T} V \rightarrow {\rm T} V [[\hbar]] \rightarrow A_{\hbar} \rightarrow A$
is the presentation ${\rm T} V \rightarrow A$. Let $I$ be the kernel of $\pi$
and $\hat{I}$ its completion in $\hat{\rm T} V[[\hbar]] = \C[[\hbar]] \,
\hattimes \, \hat{\rm T} V$. Denoting the nuclear algebra 
$\hat{\rm T} V [[\hbar]]/\hat{I} $ by $\hat{A}_{\hbar}$ we then have a 
commutative diagram
\begin{center}
\setsqparms[1`1`1`1;1000`500]
\square[ {{\rm T} V[[\hbar]]} ` {\hat{\rm T} V[[\hbar]]} ` A_{\hbar} `
{\hat{A}_{\hbar}} ; \iota` \pi `\hat{\pi}` \hat{\iota} ]
\end{center}
with injective horizontal  and surjective vertical arrows.
Only the injectivity of $\hat{\iota}$ is not immediately clear.
It follows from the fact that $I$ is closed in ${\rm T} V[[\hbar]]$:
every image of $I$ under a projection of ${\rm T} V[[\hbar]]$ to a finite 
dimensional vector space is closed. Note that the morphisms in the 
above diagram are also filtered with respect to the filtrations 
induced by the maximal ideals generated by $\hbar$. 
By these considerations we now have 
$\hat{A}_{\hbar} / \hbar \hat{A}_{\hbar} \cong\overline{A}$.
Further the vector space $A[[\hbar]] \cong A_{\hbar}$ lies densely in
$\hat{A}_{\hbar}$, hence $\hat{A}_{\hbar} \cong A[[\hbar]]$
and the algebra $\hat{A}_{\hbar}$ comprises a formal deformation of 
$\hat{A}$. 

In the next step assume $\Omega $  to be a connected open domain in $\C$
containing the origin. This gives us a continuous and injective map
$\rho: \, {\hat{\rm T} V}_{\Omega} = {\cal O} (\Omega) \hattimes 
\hat{\rm T} V \rightarrow {\hat{\rm T}}V[[\hbar]]$
by power series expansion around the origin. We now say that the formal 
deformation $A_{\hbar}$ has {\bf holomorphic initial data} if there exists a 
system of generators $Y$ of $I$ such that $Y \in {\rm im}(\rho)$.
Assuming this is the case indeed let $\tilde{Y}$ be the preimage of $Y$
under $\rho$. Let $I_{\Omega}$ be the closed ideal in 
${\hat{\rm T} V}_{\Omega}$ generated by $\tilde{Y}$. 
We then have $I_{\Omega} = \rho^{-1}( \hat{I})$,
which induces an injective and filtered homomorphism
$\tilde{A} := {\cal O} (\Omega) \hattimes \hat{\rm T} V / 
I_{\Omega} \rightarrow 
\hat{\rm T} V [[\hbar]] / \hat{I} \cong \hat{A} [[\hbar]]$.
This map is surjective from $\m^n \tilde{A} /\m^{n+1} \tilde{A}$
to $\hbar^n \hat{A} [[\hbar]] / \hbar^{n+1} \hat{A} [[\hbar]]$
where $\m$ is the maximal ideal in ${\cal O} (\Omega)$ of functions vanishing
at the origin. Thus $\tilde{A} / \m \tilde{A} \cong \hat{A}$.
As furthermore $\tilde{A}$ is dense in $\hat{A}[[\hbar]]$ we
finally have $\tilde{A} \cong \hat{A}_{\Omega}$ as a nuclear space.
Thus $\tilde{A}$ is a topologically free holomorphic deformation
of $\hat{A}_{\Omega}$.

We subsume these results in the following proposition.
\begin{Proposition}
\label{CoHoDeFoDe}
  Let $A$ be a finitely generated complex algebra (resp.~bialgebra or Hopf 
  algebra) and let $\hat{A}$ be the completion of $A$ with 
  respect to the topology of finite dimensional representations
  (resp.~nuclear Fr\'echet representations).
  Then every formal deformation  $A_{\hbar}$ of $A$ induces a formal 
  deformation $\hat{A}_{\hbar}$ of $\hat{A}$ together with a
  canonical filtered embedding $A_{\hbar} \rightarrow 
  \hat{A}_{\hbar}$. 
  If the deformation $A_{\hbar}$ has holomorphic initial data over some 
  open connected domain $\Omega \subset \C$ containing the origin, 
  then there exist a topologically free holomorphic deformation 
  $\tilde{A}$ of $A$ over $\Omega$ and a canonical filtered embedding 
  $\tilde{A}_{\Omega} \rightarrow \hat{A}_{\hbar}$. 
  These constructions are unique up to isomorphy.
\end{Proposition}
\begin{rRemark}
\label{FoToRMa}
  Suppose $H_{\hbar}$ to be a formal deformation of a finitely generated
  Hopf algebra $H$ with holomorphic initial data. 
  Give $H_{\hbar}$ the projective limit topology with respect to all
  projections $H_{\hbar} \rightarrow  H_{\hbar} / \hbar^m H_{\hbar} \cong H^m
  \subset \hat{H}^m$.
  Assume further that $H_{\hbar}$ is topologically (quasi-) triangular
  with universal R-matrix 
  ${\cal R}_{\hbar} \in H_{\hbar} \hattimes H_{\hbar}$.
  The question now is whether the universal R-matrix  on $H_{\hbar}$
  can be pushed down to one on $\tilde{H}_{\Omega}$ or in other words
  whether it induces on $\tilde{H}_{\Omega}$ the structure of a 
  topologically (quasi-) triangular Hopf algebra. 
  Obviously this is the case if and only if the formal power series
  expansion of ${\cal R}_{\hbar}$ with respect to $\hbar$
  in $\hat{H} \hattimes \hat{H}$ is converging over  $\Omega$.                      
\end{rRemark}
\begin{rRemark}
  Let $\Omega $ be an open domain in $\C^n$ or even more generally a 
  complex manifold. Then it is possible to define holomorphic deformations
  over $\Omega$ exactly like in the case of one complex variable.
  Such deformations are called {\bf holomorphic multiparameter deformations}. 
  We do not give examples here but mention that a formal multiparameter
  deformation (cf.~{\sc Reshetikhin} \cite{Res:MuQu}) with 
  holomorphic initial data obviously induces a holomorphic multiparameter 
  deformation. The theory of holomorphic deformations give sense also 
  for infinite dimensional domains $\Omega$, at least if 
  ${\cal O} (\Omega)$ is a nuclear Fr\'echet algebra. This is 
  the case, for example, if $\Omega$ is a domain in a suitable K\"othe 
  sequence space $E$ which is the dual of a Fr\'echet nuclear space.  
  (cf.~{\sc Schottenloher} \cite{Sch:MiPrAlHoFu}).
\end{rRemark}

\section{Quantum groups as holomorphically deformed Hopf algebras}
\label{QuGrHoHoAl}
By discussing some important examples we want to convince the reader 
in this section that the theory of quantum groups can be understood as the 
theory of holomorphic deformations of universal enveloping algebras
respectively of Hopf algebras of functions on a Lie group. In particular,  
we want to show that these two approaches are dual to each other in the 
context of holomorphic deformation. But this is the precise mathematical 
meaning of the physicists claim that quantizing the algebra of matrix
coefficients on a (Poisson) Lie group is dual to quantizing the corresponding 
universal enveloping algebra. 

\begin{rPar}
  {\bf Quantized ${\cal U} {\fr sl} (N+1,\C)$.}
  Consider the Lie algebra ${\fr sl} (N+1, \C)$. Its Lie algebra structure
  is given by the basis $X^+_{\iota},X^-_{\iota},H_{\iota}$ for 
  $1 \leq \iota \leq N$ together with the commutation relations
\begin{eqnarray}
\label{CoSlC1}
 & [ H_{\iota} , H_{\kappa} ] =0 \quad [H_{\iota} , X^{\pm}_{\kappa}] = \pm 
 a_{\iota \kappa} X^{\pm}_{\kappa} \\
\label{CoSlC2}
 & [X_{\iota}^+ , X_{\kappa}^- ] = \delta_{\iota \kappa} \, H_{\iota} , \quad
 [X_{\iota}^{\pm} , X_{\kappa}^{\pm} ] =0 \mbox{ for } |\iota - \kappa | >1 &\\
\label{CoSlC3}
 & {X_{\iota}^{\pm}}^2 X_{\kappa}^{\pm} - 2 X_{\iota}^{\pm} X_{\kappa}^{\pm} 
 X_{\iota}^{\pm} + X_{\kappa}^{\pm} {X_{\iota}^{\pm}}^2 = 0 \mbox{ for }
 | \iota - \kappa | = 1. &
\end{eqnarray}
  where $(\alpha_{\iota \kappa})$ is the Cartan matrix of 
  ${\fr sl} (N+1,\C)$, i.e.~$a_{\iota \iota} =2$, $a_{\iota \kappa} = 0$ 
  for $|\iota - \kappa| > 1$ and $a_{\iota \kappa} =1$ for
  $|\iota - \kappa| =1$. From the universal enveloping algebra 
  ${\cal U}{\fr sl}(N+1,\C)$ we can construct by completion with respect to
  the projective limit topology of all finite dimensional representations 
  (resp.~all nuclear Fr\'echet representations) the nuclear algebra 
  resp.~nuclear Hopf algebra $\hat{{ \cal U}}{\fr sl}(N+1,\C)$ 
  (cf.~\ref{UnEnAl}).

  We now want to construct a holomorphic deformation
  of $\hat{{\cal U}} {\fr sl}(N+1,\C)$ 
  (resp.~one of $\check{{\cal U}} {\fr sl}(N+1,\C)$) over the domain 
  $\Omega = \left\{z \in \C: \:z \neq  k \pi i, \: k \in \Z^* \right\}$
  by applying Proposition \ref{CoHoDeFoDe}. It is a well-known fact from the  
  theory of quantum groups (cf.~{\sc Chari, Pressley} \cite{ChaPre:GQG}) 
  that the following relations define a formal deformation 
  of the nuclear Hopf algebra ${\cal U} {\fr sl} (N+1, \C)$:
\begin{eqnarray}
\label{CoqSlC1}
  & [ H_{\iota} , H_{\kappa} ] =0, \quad [H_{\iota} , X^{\pm}_{\kappa}] = \pm 
  a_{\iota \kappa} X^{\pm}_{\kappa}, & \\
\label{CoqSlC2}
  & [X_{\iota}^+ , X_{\kappa}^- ] = \delta_{\iota \kappa} \, 
  \frac{\sinh \frac{z}{2} H}{\sinh \frac{z}{2}} , \quad
  [X_{\iota}^{\pm} , X_{\kappa}^{\pm} ] =0 \mbox{ for } |\iota - \kappa | >1, &\\
\label{CoqSlC3}
  & {X_{\iota}^{\pm}}^2 X_{\kappa}^{\pm} - \left( \e^{z/2} + \e^{-z/2} \right)
  X_{\iota}^{\pm} X_{\kappa}^{\pm} 
  X_{\iota}^{\pm} + X_{\kappa}^{\pm} {X_{\iota}^{\pm}}^2 = 0 \mbox{ for }
  | \iota - \kappa | = 1, & \\[2mm]
\label{CoqSlC4}
  & \Delta (H_{\iota} ) = H_{\iota} \otimes 1 + 1 \otimes H_{\iota},& \\
\label{CoqSlC5}
  & \Delta (X_{\iota}^+) = X_{\iota}^+ \otimes \e^{z \, H_{\iota}} + 1 \otimes
  X_{\iota}^+, \quad 
  \Delta (X_{\iota}^-) = X_{\iota}^- \otimes 1 + \e^{- z \, H_{\iota}}
  \otimes X_{\iota}^-, & \\[2mm]
\label{CoqSlC6}
  & S(H_{\iota}) = - H_{\iota}, \quad S(X_{\iota}^+ ) = - X_{\iota}^+ \, 
  \e^{-z\, H_{\iota}} , \quad S(X_{\iota}^-) = - \e^{z \, H_{\iota} } 
  X_{\iota}^-, & \\[2mm]
\label{CoqSlC7} 
  & \varepsilon (H_{\iota}) = \varepsilon (X_{\iota}^{\pm}) = 0.&
\end{eqnarray}
  By power series expansion around $0$ it is clear that 
  $\frac{\sinh \frac{z}{2} H}{\sinh \frac{z}{2}}$ comprises a holomorphic
  function on $\Omega$; the other relations are obviously holomorphic in $z$
  as well.
  Hence, the formal deformation of ${\cal U} {\fr sl} (N+1, \C)$ has 
  holomorphic
  initial data. By Proposition \ref{CoHoDeFoDe} the relations
  (\ref{CoqSlC1}) to (\ref{CoqSlC7})
  thus generate a holomorphic deformation 
  ${\cal U}_{q} {\fr sl} (N+1, \C)=\hat{\cal U}_{q} {\fr sl} (N+1, \C)$ 
  of the nuclear Hopf algebra $\hat{{\cal U}} {\fr sl} (N+1, \C)$
  and a holomorphic deformation 
  $\check{\cal U}_{q} {\fr sl} (N+1, \C)$ of 
  $\check{{\cal U}} {\fr sl} (N+1, \C)$.

  Using the Drinfeld double (cf.~{\sc Drinfel'd} \cite{Dri:QG})
  one can construct a (topological) R-matrix ${\cal R}_{\hbar}$ on 
  formally quantized ${\cal U} {\fr sl} (N+1, \C)$ such that
  ${\cal R}_{\hbar}$ has an expansion of the form
\begin{equation}
  {\cal R}_{\hbar} = \sum_{\beta \in \N^N} \, \left( \exp \left( 
  \hbar \left[ \frac 12 t_0  + \frac 14 (H_{\beta} \otimes 1 + 
  1 \otimes H_{\beta} )\right] \right) \right) P_{\beta}, 
\end{equation}
  where $t_0 \in {\fr sl} (N+1, \C) \otimes {\fr sl} (N+1, \C)$ is 
  chosen appropriate, $H_{\beta} = \sum_{\iota} \, \beta_{\iota} H_{\iota}$ 
  and the $P_{\beta}$ are polynomials homogeneous of degree 
  $\beta_{\iota}$ in  $ X^+_{\iota} \otimes 1$ and  $1 \otimes X^-_{\iota}$.   
  Hence ${\cal R}_{\hbar}$ has a converging power series expansion, so
  Remark \ref{FoToRMa} entails that ${\cal U}_{q} {\fr sl} (N+1, \C)$
  is topologically quasitriangular as well.
\end{rPar}
\begin{rRemark}
  The quantum group models defined by {\sc Jimbo} \cite{Jim:QAUGLHAYBE} do 
  not comprise a holomorphic deformation of 
  $\hat{\cal U} {\fr sl} (N +1, \C)$,
  though a holomorphic deformation of a certain extended Hopf algebra,
  for example for $N = 1$, of  
  $\hat{\cal U} {\fr sl} (2, \C) \otimes \C[T]/(T^2-1)$ 
  (\cite{BonFlaGerPin:HGSQGSDRPD}).
\end{rRemark}
\begin{rPar}
  {\bf Quantizing $ {\rm SL} (N,\C)$ by dualizing.}
  Now we want to quantize the dual of the preceding example, or in other
  words the algebra ${\cal R}( {\rm SL} (N,\C))$ of representation functions 
  on the Lie group $ {\rm SL} (N,\C)$. Recall from \ref{SiLiAl} that
  $\hat{\cal U} {\fr sl} (N, \C)$ and $ {\cal R}( {\rm SL} (N,\C))$ 
  are strictly nuclear and that
  $\hat{\cal U} {\fr sl} (N, \C)$ is topologically isomorphic to 
  $ {\cal R}( {\rm SL} (N,\C))'$. Hence, according to
  Proposition \ref{CoHoDeFoDe} the deformation 
  ${\cal U}_{q} {\fr sl} (N, \C)$ gives rise to a holomorphic 
  deformation ${\cal U}_{q} {\fr sl} (N, \C)'$
  of  $ {\cal R}( {\rm SL} (N,\C))$. Unfortunately we do not yet have 
  a concrete representation of ${\cal U}_{q} {\fr sl} (N, \C)'$
  by generators and relations. As will be shown in the next paragraph,
  the well-known FRT-deformation of ${\cal R}({\rm SL} (N,\C))$
  does exactly provide that. 
\end{rPar}

\begin{rPar} 
 {\bf Quantizing $ {\rm SL} (N,\C)$ according to FRT.}
  First let us recall the FRT-construction of quantized 
  ${\cal R}( {\rm SL} (N,\C))$.  Any finite dimensional representation
  of ${\rm SL} (N,\C)$ can be realized as a subrepresentation of a sum
  of tensor products of the fundamental representation of ${\rm SL} (N,\C)$.
  This means that ${\cal R}( {\rm SL} (N,\C))$ is generated by the matrix
  coefficients $T = \left( t_i^j \right)_{1\leq i,j \leq N}$ of the 
  fundamental representation. The $t_i^j$ hereby 
  fulfill the relation
\begin{equation}
  {\rm det} \, T \: = \: 1.
\end{equation}
  It can even be shown that the algebra generated by $t_i^j$ modulo
  the determinant relation is in fact isomorphic to 
  ${\cal R}( {\rm SL} (N,\C))$. From \ref{MaCoGrRel} and \ref{SiLiAl} we infer 
  the finest locally convex topology on ${\cal R}( {\rm SL} (N,\C))$
  being the natural one, in particular, because of the natural duality 
  between ${\cal U} \g$ and ${\cal R} (G)$. To deform the nuclear Hopf algebra
  ${\cal R}( {\rm SL} (N,\C))$ 
  we will now use the FRT-construction as described in \ref{FRT}. 
  Let $V$ be the complex vector space spanned by $x_i$, $i = 1,...,N$ and
  consider the $R$-matrix 

\begin{eqnarray}
\label{SLNRMa}
  R : \: \C^* & \rightarrow &\End (V)' \otimes \End (V)', \nonumber \\ 
  z & \mapsto & 
  z \, \sum \limits_{i=1}^n \,  t_i^i \otimes t_i^i \, + \,
  \sum \limits_{i,j=1 \atop i \neq j}^n \,  t_i^i \otimes t_j^j \, + \,
  (z - z^{-1}) \,  \sum \limits_{i,j=1 \atop i > j}^n \,  
  t_i^j \otimes t_j^i,
\end{eqnarray} 
  where $(t_i^j)$ is the basis of ${\rm End} \, (V)'$ dual to the basis 
  $(e_i^j)$ of ${\rm End} \, (V)$, and $e_i^j$ is the endomorphism of $V$ 
  mapping $x_i$ to $x_j$ and
  vanishing on all the other elements of $(x_1,...,x_N)$.
  This $R$-matrix fulfills the quantum Yang-Baxter equation and is
  nondegenerate with inverse $R(z) ^{-1} = R (z^{-1})$. 
  (see for example {\sc Takhtajan} \cite{Tak:LQG}).
  Hence, by the FRT-construction \ref{FRT} we receive the bialgebra $A(R)$.
  It is well-known that the {\bf quantum determinant}
\begin{equation}
  {\rm det_q} \, T \: = \: \sum \limits_{\sigma \in {\rm S}_N} \,
  (- z)^{\ell (\sigma)} \, t_1^{\sigma (1)} \cdot ... \cdot t_N^{\sigma (N)}   
\end{equation}   
  belongs to the center of $A(R)$. Denoting by $I$ the ideal 
  $I = A(R) ({\rm det_q} \, T -1 ) $ the quotient bialgebra 
  $A (R) / I $ then is even a Hopf algebra with antipode given by
\begin{equation}
  S (t_i^j) \: = \: (-z)^{i-j} \, \tilde{t}_i^j,
\end{equation}
  where the $\tilde{t}_i^j$ are the so-called quantum-cofactors
\begin{equation}
  \tilde{t}_i^j \:  = \:  \sum \limits_{\sigma \in {\rm S}_{N-1}} \,
  (- z)^{\ell (\sigma)} \, t_1^{\sigma_{1}} \cdot ...\cdot
  t_{i-1}^{\sigma_{i-i}} \cdot t_{i+1}^{\sigma_{i+i}} \cdot ...
  \cdot t_N^{\sigma_{N}}  
\end{equation}
  with
\begin{equation}
  (\sigma_1, ..., \sigma_{i-1}, \sigma_{i+1} ,..., \sigma_N) \:  =
  \:  (\sigma (1), ..., \sigma (j-1), \sigma (j+1) ,..., \sigma (N)). 
\end{equation}
  We denote the quotient $A (R) / I $ by 
  ${\cal R} ({\rm SL_q } (N , \C))$ and call it the 
  {\bf algebra of matrix coefficients on quantized ${\rm SL} (N,\C)$}.
  As we will show ${\cal R} ({\rm SL_q } (N , \C))$ (endowed with the 
  locally convex inductive limit topology of all finite dimensional subspaces)
  comprises the holomorphic deformation of 
  ${\cal R} ({\rm SL} (N , \C))$ we are looking for. In particular 
  ${\cal R} ({\rm SL_q } (N , \C))$ will be isomorphic to
  ${\cal U}_{q}' {\fr sl} (N, \C)$.
  \vspace{2mm}

  In the second step we will show that the FRT-bialgebra $A(R)$ 
  is a holomorphic deformation of the bialgebra of polynomial functions 
  on the semigroup ${\rm M} (N \times N, \C)$. 
  Though the proof seems to be quite technical we will give it in some detail,
  as it provides an example for constructing a holomorphic deformation
  without first having a formal one.   
  
  Let us determine the relations between the $t_i^j$. These
  relations are given by the entries of the matrix 
  $\: T_1 \, T_2 \, R  - R \, T_2 \, T_1$,  where $T_1 = T \odot 1$ and 
  $T_2 = 1 \odot T$. 
  Equivalently we can calculate the entries of the matrix
\begin{equation}
  T_2 \, T_1 \: - \: R(z^{-1}) \, T_1 \, T_2 \, R (z).
\end{equation}
  They are written down in the following table:
\begin{eqnarray}
\label{CoReSLBiAl}
\begin{array}{|lc|cl|}
  \hline  & & & \\[-3mm] 
  \mbox{variables } i,j,k,l & & & t_i^j \otimes t_k^l \: - \: 
  \big( R(z^{-1}) \, T_1 \, T_2 \, R (z) \big)_{ik,jl} \\[2mm] 
  \hline & &   & \\[-3mm] 
  i = k, \: j < l & & & t_i^j \otimes t_i^l \: - \: z \, t_i^l \otimes t_i^j \\[2mm]
  i = k, \: j > l & & & t_i^j \otimes t_i^l \: - \: z^{-1} \, t_i^l \otimes t_i^j \\[2mm]
  i < k, \: j =l & & & t_i^j \otimes t_k^j \: - \: z \, t_k^j \otimes t_i^j \\[2mm]
  i > k, \: j =l & & & t_i^j \otimes t_k^j \: - \: z^{-1} \, t_k^j \otimes t_i^j \\[2mm]
  i < k, \: j>l & & & t_i^j \otimes t_k^l \: - \: t_k^l \otimes t_i^j  \\[2mm]
  i < k, \: j<l & & & t_i^j \otimes t_k^l \: - \: \left( t_k^l \otimes t_i^j 
                      + (z - z^{-1}) \, t_k^j \otimes t_i^l \right) \\[2mm]
  i > k, \: j > l & & & t_i^j \otimes t_k^l  \: - \: \left( t_k^l \otimes t_i^j
                      - (z - z^{-1}) \, t_i^l \otimes t_k^j \right) \\[2mm]
  i > k, \: j < l & & & t_i^j \otimes t_k^l \: - \\ & & &
          \hspace{2mm} \big( t_k^l \otimes t_i^j + 
          (z - z^{-1}) \left( t_k^j \otimes  t_i^l  - t_i^l \otimes  t_k^j - 
          (z - z^{-1} ) \, t_i^j \otimes t_k^l \right) \big) \\[2mm] \hline
\end{array}
\end{eqnarray} 
Check that the first, the third, the sixth and last entry in this table are 
linear combination of the other ones.
Hence, $A(R)$ is generated by the $t_i^j$ modulo the relations
\begin{eqnarray}
\left\{ \begin{array}{l}
\left( t_i^j \otimes t_k^l \: - \: z^{-1} \, t_k^l \otimes t_i^j 
       \right)_{i = k, \: j > l} , \:
\left( t_i^j \otimes t_k^l \: - \:  t_k^l \otimes t_i^j 
       \right)_{i > k, \: j < l }, \\ 
\left( t_i^j \otimes t_k^l \: - \:\left(  t_k^l \otimes t_i^j 
       - (z - z^{-1}) \, t_k^j \otimes t_i^l \right) 
       \right)_{ i > k, \: j > l} ,  
\left( t_i^j \otimes t_k^l \: - \: z^{-1} \, t_k^l \otimes t_i^j  
       \right)_{ i > k,  \: j = l}
\end{array} 
\right\}.
\end{eqnarray}  
Now it is straightforward to show that $A(R)$ is free over
${\cal O} (\C^*)$ with basis 
\begin{equation}
\label{BaSLq}
\big( T^m \big)_{m \in \N^{N^2}} \: := \: 
 \left( \big( t_1^1 \big)^{m_1^1} \cdot
 ... \cdot \left( t_1^N \right)^{m_1^N} \cdot \big( t_2^1 \big)^{m_2^1} 
 \cdot ... \cdot  \left( t_N^N \right)^{m_N^N}
 \right)_{m = (m_1^1,m_1^2,...,m_N^N) \in \N^{N^2} } .
\end{equation}
Consequently $A(R)$ is isomorphic as ${\cal O} (\C^*)$-module  to 
$ {\cal O} (\C^*) \otimes A(R)_1 = {\cal O} (\C^*) \hattimes A(R)_1 $,
where the undeformed algebra $A(R)_1 = A(R) / \m_1 A(R)$ is 
isomorphic to the algebra $\C [ t_1^1 , ..., t_1^N,..., t_N^N]$
of polynomial functions on ${\rm M} (N \times N , \C)$. The $t_i^j$ 
can hereby be interpreted as the functions giving the $(i,j)$ entry of a
matrix. Altogether these considerations prove that $A(R)$ is a holomorphic
deformation of $\C [ t_1^1 , ..., t_1^N,..., t_N^N]$ indeed.
Note that $A(R)$ does not have zero divisors.
\begin{rRemark}
  The basis of $A(R)$ given in (\ref{BaSLq}) could as well be derived
  using the diamond lemma of {\sc Bergman} \cite{Ber:DLRT}.
\end{rRemark}

Next let us prove that the monomials $  \big( T^m + I \big)_{m \in M}$  with 
\begin{equation}
 M = \left\{ (m_1^1,..., m_N^N) \in \N^{N^2} : \mbox{ one element of }
 \left\{ m_1^N,...,m_N^1 \right\} \mbox{ vanishes} \right\} 
\end{equation}
form a basis of ${\cal R} \left( {\rm SL_q} (N ,\C) \right) = A(R) / I$.
Denote for $m \in \N^{N^2}$ by ${\rm val }\, m$ the {\bf value}
\begin{equation}
  {\rm val} \, m \: = \: m_1^N + ... + m_N^1.
\end{equation}
To prove that ${\rm span} \left\{ \left( T^m + I \right)_{m \in M}\right\}
= A(R) /I$ it suffices  to show that for $m \in \N^{N^2} \setminus M$ we can 
write 
the  monomial $T^m $  modulo $I$ as a ${\cal O} (\C^*)$-linear combination of 
monomials $T^r $, $r \in \N^{N^2}$ with ${\rm val} \, r < {\rm val} \, m$.
Indeed, repeating this process for the $r$ and proceeding inductively
we can after finitely many steps expand $T^m $ modulo $I$
as a combination of monomials $T^s $, $s \in \N^{N^2}$
with ${\rm val } \, s < {\rm val } \, m$ and $ s \in M$.
The idea now is to move for any pair $(i,j)$ with $i+j = N+1$ one $t_i^j$
appearing in $T^m $ to the right while controlling the value of 
the monomials created by this process. Denote by $a_k^l \in \N^{N^2}$  
the multiindex the entry of which at the $\left( _k^l \right)$ position
is $1$ and $0$ otherwise. Suppose we already have for an index
$i \in \{ 1,...,N\}$ the expansion
\begin{equation}
\label{ExMo}
  T^m  \: = \: \sum \limits_{n \in \N^{N^2} } \, c(n,i) \, T^n \, + \,
  d(i) \,  T^{m(i)} \cdot t_{i}^{N+1-i} \cdot ... \cdot t_N^1,
\end{equation}
where $c(n,i) , d(i) \in {\cal O} (\C^* )$, ${\rm val} \,n <  {\rm val} \,m$
for every $n \in \N^{N^2}$ with $c(n,i) \neq 0$
and $m(i) = m - \sum_{k \geq i} \, a_k^{N+1-k}$.
Now consider one $t_{i-1}^{N+2-i}$ in $T^m$. By the commutation relations 
(\ref{CoReSLBiAl}) we have
\begin{equation}
  T^{m(i)} \: = \: \tilde{d} (i) \, T^{m(i) - a_{i-1}^{N+2-i}} \, 
  t_{i-1}^{N+2-i} \, + \, \sum \limits_{n \in M(i)} \, \tilde{c} (n,i) \, T^n,
\end{equation} 
  where $\tilde{d} (i), \tilde{c} (i,n) \in {\cal O} (\C^*)$ and
\begin{equation}
  M(i) \: = \: \left\{ n \in \N^{N^2} : 
\begin{array}{l}
  n_1^1 + n_1^2 + ... + n_N^N =
  m(i)_1^1 + m(i)_1^2 + ... + m(i)_N^N,  \\
  n_i^{N+1-i} = m(i)_i^{N+1-i} -1 , \quad n_k^{N+1-k} = m(i)_k^{N+1-k}
  \mbox{ for } k \neq i
\end{array}
\right\}.
\end{equation}
By definition ${\rm val }\, n < {\rm val }\, m(i)$ holds for every 
$n \in M(i)$. Hence, by (\ref{CoReSLBiAl}) the expansion of 
$T^n \cdot  t_{i-1}^{N-i} \cdot ... \cdot t_1^N$, $n \in M(i)$ 
with respect to the basis $\big( T^{\tilde{m}} \big)_{\tilde{m} \in \N^{N^2}}$
contains only monomials $T^{\tilde{m}}$ with 
${\rm val} \, \tilde{m} < {\rm val} \, m $.
Therefore the relation 
\begin{equation}
  T^m  \: = \: \sum \limits_{n \in \N^{N^2} } \, c(n,i-1) \, T^n \, + \,
  d(i-1) \,  T^{m(i-1)} \, t_{i-1}^{N+2-i} \cdot ... \cdot t_N^1
\end{equation}
is true, and ${\rm val} \, n < {\rm val} \, m $  for every $n \in \N^{N^2}$ 
with $c(n,i-1) \neq 0$. 
Hence, by induction we can set $i = 1$ in (\ref{ExMo}).
Modulo $I$ this implies that
\begin{eqnarray}
  T^m  & = & \sum \limits_{n \in \N^{N^2} } \, c(n,1) \, T^n \, + 
  \, d \,  T^{m(1)} \, + \,
  \sum \limits_{\sigma \in S_n \setminus \{ \sigma_{\rm inv} \}}
  d(\sigma) \, T^{m(1)} \, 
  t_{1}^{\sigma (1)} \cdot ... \cdot t_N^{\sigma (N)},
\end{eqnarray}
where $\sigma_{\rm inv} $ is the permutation fulfilling
$\sigma_{\rm inv} (i) = N+1-i$ for every $i=1,...,N$ and
$c(n), d, d(\sigma) $ are elements of ${\cal O} (\C^*)$ such that
$c(n) = 0$ for $n \in \N^{n^2}$ with ${\rm val} \, n  \geq {\rm val }  \, m$.
Now we want to determine the expansion of 
$T^{m(1)} \,  t_{1}^{\sigma (1)} \cdot ... \cdot t_N^{\sigma (N)}$
with respect to the basis $\left( T^{\tilde{m}} \right)_{\tilde{m} \in \N^{N^2}}$
by ``moving the $t_i^{\sigma (i)}$ to the left''. 
Check that for every $\sigma \in S_n \setminus \{ \sigma_{inv} \}$ there exists
an index $i$ such that $i + \sigma (i) > N+1$.  
Moving the corresponding $t_i^{\sigma (i)}$ in the monomial
$ T^{m(1)} \,  t_{1}^{\sigma (1)} \cdot ... \cdot t_N^{\sigma (N)} $
to the left by using the commutation relations (\ref{CoReSLBiAl})
one does not generate a further $t_k^{N+1-k}$. 
Moving another $t_j^{\sigma (j)}$ to the left one generates
at most one $t_k^{N+1-k}$. This entails that 
$ T^{m(1)} \,  t_{1}^{\sigma (1)} \cdot ... \cdot t_N^{\sigma (N)} $
can be expanded by monomials $ T^{\tilde{m}}$ with
${\rm val} \, \tilde{m} \leq m(1) + N -1 = {\rm val} \, m - 1$.
But then $T^m$ itself can modulo $I$ be expanded by such monomials,
hence ${\rm span} \left\{ \left( T^m + I \right)_{m \in M}\right\} = A(R) /I$.

We still have to prove that $\left( T^m + I \right)_{m \in M}$ is a linear
independant family. So let
\begin{equation}
\label{CoEqExMo}
  \sum \limits_{m \in M} \, d(m) \, T^m \:  = \:
  a \, \left( {\rm det_q} \, T - 1 \right)
\end{equation} 
with $d(m) \in {\cal O} (\C^*) $ and $a \in A(R)$,
and assume that both sides do not vanish.
But then by (\ref{CoReSLBiAl}) and the definition of ${\rm det_q}T$ 
the expansion of the right hand side contains at least
one monomial $T^m$, $m \notin M$ with nonvanishing coefficient in 
${\cal O}(\C^*)$.
Thus the left and right side of (\ref{CoEqExMo}) cannot be equal
unlike they are both $0$. This shows the linear independance.

As the monomials $T^m$, $m \in M$ form a $\C$-basis of the 
undeformed algebra ${\cal R} ({\rm SL} (N,\C))$, it follows that 
as (topological) vector spaces  
${\cal O} (\C^*) \otimes {\cal R} ({\rm SL} (N,\C))$ and 
${\cal R} ({\rm SL_q} (N,\C))$ are isomorphic.
Furthermore we have the following chain of Hopf algebra isomorphisms:
\begin{eqnarray}
\lefteqn{
 {\cal R} ({\rm SL_q} (N,\C)) / \m_1 {\cal R} ({\rm SL_q} (N,\C)) \cong }
 \nonumber \\
 & \cong & A(R) / \m_1 A(R) \, \big/ \, I / \m_1 I 
 \cong  A(R)_1 \, / \,A(R)_1  ({\rm det } \, T -1 )
 \cong  {\cal R} ({\rm SL} (N,\C)).
\end{eqnarray}
This finally proves the claim, i.e.~${\cal R} ({\rm SL_q} (N,\C)) $
is a holomorphic deformation of ${\cal R} ({\rm SL} (N,\C)) $. \vspace{2mm}
\vspace{2mm}

Next we want to sketch the proof for the isomorphy of 
${\cal R}  ({\rm SL_q} (N,\C)) $ and ${\cal U}_q {\fr sl} (N, \C)'$
or in other words for the duality between
${\cal R}  ({\rm SL_q} (N,\C)) $ and ${\cal U}_q {\fr sl} (N, \C)$.
Denote for $z \in \C$ by ${\cal R}  ({\rm SL}_z (N,\C)) $ the Hopf algebra 
${\cal R}  ({\rm SL_q} (N,\C)) \big/ \m_z \, {\cal R}  ({\rm SL_q} (N,\C)) $.
Similarly let ${\cal U}_z {\fr sl} (N, \C) = {\cal U}_q {\fr sl} (N, \C)
\big/ \m_z \, {\cal U}_q {\fr sl} (N, \C)$ for $z \in \Omega$.
We claim that ${\cal R}  ({\rm SL}_{{\rm e}^z} (N,\C))$, $z \in \Omega$
is topologically isomorphic to the restricted Hopf dual 
${\cal U}_z {\fr sl} (N, \C)^{\circ}$ with the finest locally convex
topology. Check that there is a canonical $N$-dimensional indecomposable 
${\cal U}_z {\fr sl} (N, \C)$-module 
$\rho: {\cal U}_z {\fr sl} (N, \C) \rightarrow \End (W)$.
Define the matrix coefficients $A_{ij}$ by
\begin{eqnarray}
 \rho (h) & = & \left(
\begin{array}{ccccc}
  A_{11} (h) & \cdot & \cdot & \cdot & A_{1n} (h) \\
  \cdot & \cdot & & & \cdot \\
  \cdot & & \cdot & & \cdot \\
  \cdot & & & \cdot & \cdot \\
  A_{n1}(h) & \cdot & \cdot & \cdot & A_{nn} (h)
\end{array} \right)
\end{eqnarray}
for any $h \in {\cal U}_z$. Therefore, we can define a homomorphism
$\Phi: {\rm T} C \rightarrow {\cal U}_z {\fr sl} (N, \C)^{\circ}$,
where $C = {\rm End} (W)'$. Obviously, $\Phi$ is continuous with respect
to the finest locally convex topology on ${\rm T} C$.
By a lengthy calculation one can show that $\Phi$ is surjective
and $\ker \Phi = I_z := \ker ( {\rm T} C \rightarrow {\cal R}
({\rm SL}_{{\rm e}^z} (N ,\C) )$.
(More details for this are given in 
{\sc Chari, Pressley} \cite{ChaPre:GQG}, Theorem 7.1.4. 
for the case of  ${\cal U} {\fr sl} (2, \C)$ and
in {\sc Takeuchi} \cite{Tak:STGL,Tak:HATAQGUSL} for the general case.
See also {\sc Drinfel'd} \cite{Dri:QG}.)
Hence, $\Phi$ induces a linear isomorphism 
$\overline{\Phi}: {\cal R} ({\rm SL}_{{\rm e}^z} (N, \C))  \rightarrow 
{\cal U}_z {\fr sl} (N, \C)^{\circ}$ which is topological since both spaces
are endowed with the finest locally convex topology.

We now have proven the main part of the following theorem.
\begin{Theorem}
  The FRT-algebra 
  ${\cal R} ({\rm SL_q} (N, \C))  = A(R)/(\det_{\rm q} -1 ) A(R)$
  corresponding to the $R$-matrix (\ref{SLNRMa})
  comprises a holomorphic quantization of the Poisson algebra 
  ${\cal R} ({\rm SL} (N, \C)) $ of matrix coefficients on the Lie group
  ${\rm SL} (N,\C)$. Moreover, it coincides 
  with the deformation ${\cal U}_{\rm q} {\fr sl} (N, \C)'$ 
  dual to the quantization of ${\cal U}{\fr sl} (N, \C)$. 
  For every $z \in \Omega$ the Hopf algebra 
  ${\cal R} ({\rm SL}_{{\rm e}^z} (N, \C))$ is topologically isomorphic 
  to the restricted Hopf dual ${\cal U}_z {\fr sl} (N, \C)^{\circ}$, and 
  ${\cal U}_z {\fr sl} (N, \C)$ to 
  ${\cal R} ({\rm SL}_{{\rm e}^z} (N, \C))'$.
\end{Theorem}
\begin{Proof}
  It has been shown above that $A(R) / (\det_{\rm q} -1 ) A(R)$
  is a holomorphic deformation of ${\cal R} ({\rm SL} (N, \C)) $.
  That it is even a holomorphic quantization follows from
  Theorem \ref{HQFFQ} and \cite{FadResTak:QLGLA}.
  
  As $\overline{\Phi}: 
  {\cal R} ({\rm SL}_{{\rm e}^z} (N, \C))  \rightarrow 
  {\cal U}_z {\fr sl} (N, \C)^{\circ}$ is a  topological isomorphism,
  paragraph \ref{ToSeSiUnEnAl} entails that 
  ${\cal R} ({\rm SL}_{{\rm e}^z} (N, \C))$ is isomorphic to
  $ {\cal U}_z {\fr sl} (N, \C)'$. By strict nuclearity 
  ${\cal U}_z {\fr sl} (N, \C) \cong
  {\cal R} ({\rm SL}_{{\rm e}^z} (N, \C))'$ follows as well.
\end{Proof}
\end{rPar}

\appendix
\section{Nuclear spaces and holomorphic vector-valued functions}
\label{NuSp}
\subsection*{Locally convex spaces}
A seminorm on a (real or complex) vector space $E$ is a function 
$p: E \rightarrow \R$ 
satisfying $p \geq 0$, $p (0) = 0$, $p(x+y) \leq p(x) + p(y)$ and 
$p(\lambda x) = |\lambda | p(x)$ for all $x,y \in E$ and $\lambda \in \K$. 
Any family $P$ of seminorms on $E$ defines a locally convex topology on $E$: A 
subbasis of this topology is, for example, the system of all ``open balls''
${\rm B}_p (a, r) =  \left\{ x \in E : \: p (a -x) < r \right\}$, where
$a \in E$, $r > 0$ and $ p \in P$. 
$E$ together with such a topology is called a locally convex space, and the set
of all continuous seminorms on $E$ 
will be denoted by ${\rm cs } (E)$. A locally
convex space can also be described as a vector space $E$ together with
a translation invariant topology for which the structure maps 
$E \times E \rightarrow E$ and $\K \times E \rightarrow E$ are continuous such 
that $ 0 \in E$ has a neighborhood 
basis ${\cal B}$  consisting of absolutely convex open subsets $ U \subset E$.
Each such $U$ corresponds to a seminorm 
\begin{equation}
  p_U (x) \: := \: \inf \left\{ \lambda \in \R : \: x \in \lambda U \right\},
\end{equation}
and the topology is defined by the family 
$\left\{ p_U : \: U \in {\cal B} \right\}$ of seminorms. 
In general, $P$ as above is called 
a defining family of seminorms of the locally convex space $E$.

A locally convex space $E$ is metrizable (resp.~normalizable), if $E$ is 
Hausdorff and possesses  a countable (resp.~finite) defining family  $P$ of
seminorms. $E$ is Hausdorff, if to every $x \in E \setminus \{ 0 \}$ there 
corresponds a $p \in P$ with $p (x ) \neq 0$. 
We are mainly interested in locally convex
Hausdorff spaces which are complete with respect to the uniformity given
by ${\rm cs} (E)$. A complete metrizable locally convex space is called a
Fr\'echet space.
\subsection*{Locally m-convex spaces}
A locally convex algebra over $\K$ is an associative $\K$-algebra $A$ with
a locally convex topology such that the multiplication 
$A \times A \rightarrow A$ is continuous. 
In order to have a reasonable functional calculus on $A$ one needs
an additional property, namely m-convexity. Recall that a seminorm $p$ on an
algebra $A$ is called multiplicative, if
\begin{equation}
  p (xy) \leq p(x) \, p(y)
\end{equation}
for all $x,y \in A$. A locally convex algebra $A$ is called locally m-convex
(cf.~{\sc Michael} \cite{Mic:LMCTA}) if there exists a defining family $P$ of
multiplicative seminorms. 
On a complete locally m-convex algebra over $\C$ there exists a 
functional calculus as an action of the space ${\cal O} (\C)$ of entire
functions on the algebra $A$: To each $f \in {\cal O} (\C)$ with power series
expansion $f(z) = \sum \, c_n z^n $ and each $a \in A$
the series $\sum \, c_n a^n$ yields a well-defined
element $\hat{f} (a) = \sum \, c_n a^n \in A$. The
convergence of this series follows from $p( a^n) \leq p(a)^n$ for
multiplicative seminorms $p$ on $A$: 
\begin{equation}
p \left( \sum \limits_{m \leq n \leq m + k} \right) \leq 
\sum \limits_{m \leq n \leq m + k} \, |c_n | \, p( a )^n \longrightarrow 0
\end{equation}
for $m,k \rightarrow \infty$ and $p \in P$. Hence, $ \sum \, c_n a^n $ is a
Cauchy series. An example for a locally convex but not locally m-convex
algebra is given by $\C [T]$ (or ${\cal U} \g$) with the finest locally 
convex topology. Otherwise, $\sum \, \frac{1}{n!} a^n$
would converge for all $a \in A$. But for $a = T$ the sequence 
$\sum \, \frac{1}{n!} T^n$ does not have a limit in $\C [T]$.
The usual topologies on function algebras are locally m-convex since they 
are defined by seminorms given by the 
supremum of the function (resp.~their derivatives) on a family of subsets of
the common domain of the functions in question. 
Clearly, those seminorms are multiplicative. In particular
the algebra ${\cal O} (\Omega)$ of holomorphic functions on an open
domain $\Omega \subset \C$ with the compact open topology is locally m-convex.
This holds true for an open subset $\Omega \subset \C^n$ as well, 
for a complex manifold $\Omega$ or even for an infinite dimensional domain 
in a locally convex space $E$.

Every normed algebra over $\K$ is locally m-convex, and so is every subalgebra
of a product of normed algebras. Hence, the locally convex inverse limit 
topology on $A$ of a family $\varphi_i : A \rightarrow B_i$ of homomorphisms
into normed algebras is locally m-convex, as is the locally convex projective
limit of a projective system of m-convex algebras. Hence, the projective
limit topologies considered in section \ref{NLCNHA} define locally m-convex 
algebras.
\subsection*{Nuclear Spaces}
On a given tensor product $E \otimes F$ of two locally convex spaces $E$ and 
$F$ one can consider many different locally convex topologies arising from
the topologies on $E$ and $F$. The most natural one is the $\pi$-topology, 
i.e.~the finest locally convex topology on $E \otimes F$ for which the 
natural mapping
\begin{equation}
  \otimes : \: E \times F \rightarrow E \otimes F 
\end{equation}
is continuous. $E \otimes F$ with this topology is denoted by 
$E \otimes_{\pi} F$, its completion by $E \hattimes F$.
Another useful topology can be described
by the embedding $E \otimes F \rightarrow {\cal B} (E'_{\rm s} , F'_{\rm s})$,
where $E'_{\rm s}$ is the dual $E'$ equipped with the topology of simple 
convergence. The $\varepsilon$-topology is the topology induced from 
${\cal B}_{\varepsilon} ( E'_{\rm s} , F'_{\rm s} )$, where the subscript
$\varepsilon$ denotes the topology of uniform convergence on all products
$A \times B \subset E' \times F'$ of equicontinuous subsets 
$A \subset E', B \subset F'$. Both, the $\pi$-topology and the 
$\varepsilon$-topology are compatible with $\otimes$ in the following sense:
\begin{enumerate}
\item
  $\otimes : \: E \times F \rightarrow E \otimes F$ is continuous,
\item
  for all $(e,f) \in E' \times F'$ the linear form 
  $ e \otimes f : \: E \otimes F \rightarrow \K$, $ x \otimes y \mapsto
  e(x) \, f(y)$ is continuous. 
\end{enumerate}
In fact, the $\pi$-topology is the strongest
and the $\varepsilon$-topology is the weakest topology on $E \otimes F$
compatible with $\otimes$. 

A locally convex space $E$ is called {\bf nuclear}, if all the compatible
topologies on $E \otimes F$ agree for all locally convex spaces $F$
(cf.~{\sc Grothendieck} \cite{Gro:PTTEN}). 
All finite dimensional vector spaces are nuclear. Subspaces and Hausdorff 
quotients of nuclear spaces are nuclear, as well as the completion of a 
nuclear space. 
In general, locally convex projective limits of nuclear spaces and
countable locally convex inductive limits of nuclear spaces are also nuclear.
The direct sum $\C^{(\Lambda)}$ for an uncountable $\Lambda$ provides an
example of an inductive limit of nuclear spaces which is not nuclear. It 
also provides an example of a nuclear space $E$, namely $E = \C^{\Lambda}$,
such that the dual $E ' \cong \C^{(\Lambda)}$ is not nuclear with respect
to the strong topology on $E'$, i.e.~the topology of uniform convergence on
bounded subsets of $E$.
In the following we will denote by $E \hattimes F$ the completion of
$E \otimes_{\pi} F$. Note, that $E \hattimes F$ is nuclear for 
complete nuclear spaces $E$ and $F$.  

\begin{Definition}
  We call a reflexive nuclear space $E$ {\bf strictly nuclear} if the dual
  $E'$ with the strong topology is nuclear as well and if the 
  inclusion $E' \otimes F' \rightarrow ( E \hattimes F )' $ induces an
  isomorphism 
\begin{equation}
\label{IsStNuSp}
  E' \hattimes E' \: \cong \: ( E \hattimes E) ',
\end{equation}
  algebraically and topologically. Hereby $E'$, $F'$ and $(E \hattimes F)'$
  are endowed with the strong topology which coincides with the 
  compact open topology for complete nuclear spaces $E$.  
\end{Definition}

  Formula (\ref{IsStNuSp}) holds true for nuclear Fr\'echet spaces 
  (cf.{\sc Tr\`eves} \cite{Tre:TVSDK}) and therefore, by dualizing, also for
  duals of nuclear Fr\'echet spaces. It also can be proven for nuclear 
  LF-spaces $E$ (and their duals), i.e.~nuclear spaces $E$ which can be
  described as an inductive limit $E_n \subset E_{n+1} \subset ... \subset E$
  of Fr\'echet spaces such that $E_n$ has the topology induced from $E_{n+1}$.
  As a consequence many of the important function spaces in Analysis are 
  strictly nuclear, e.g.~the space ${\cal O} (\Omega )$ of holomorphic
  functions on a complex manifold $\Omega$, the space 
  ${\cal E} (M) ={\cal C}^{\infty} (M)$ of infinitely differentiable functions
  on a manifold $M$, the space ${\cal D} (M)$ of test functions on a manifold
  $M$ and the dual, the space of distributions. The latter two examples
  are in general neither Fr\'echet nor dual Fr\'echet. Note that for a 
  strictly nuclear space $E$ the dual $E'$ is strictly nuclear as well. 
  This can be seen by simply dualizing (\ref{IsStNuSp}).

  The space $\K^{(\Lambda)}$ of functions $\Lambda \to \K$ with finite support
  for an arbitrary set $\Lambda$ is the prototype of a locally convex 
  space $E$ with the finest  locally convex topology (simply choose a Hamel 
  basis of $E$ indexed by $\Lambda$). The space $\K^{(\Lambda)}$ 
  satisfies the duality condition (\ref{IsStNuSp}) as well, although this 
  space is only nuclear (and strictly nuclear) 
  for countable $\Lambda$. The reason for this is essentially the 
  fact that for $E := \K^{(\Lambda)}$ the tensor product $E \otimes E$ 
  with the $\pi$-topology is already complete and isomorphic to 
  $\K^{(\Lambda \times \Lambda)}$, which follows from the observation 
  that the product topology on $E \times E$ is the finest locally convex 
  topology. Moreover, the strong dual 
  of $\K^{(\Lambda)}$ is isomorphic to $\K^{\Lambda}$ with 
  the product topology. Hence, 
  $$ (E \hat{\otimes} E)' \cong (\K^{(\Lambda \times \Lambda)})' \cong
    \K^{\Lambda \times \Lambda} \cong 
    \K^{\Lambda} \hat{\otimes} \K^{\Lambda} \cong
    (\K^{(\Lambda)})' \hat{\otimes} (\K^{(\Lambda)})',$$
  i.e.
\begin{equation}
\label{IsKL}
  (\K^{(\Lambda)})' \hat{\otimes} (\K^{(\Lambda)})' \cong 
  (\K^{(\Lambda)} \hattimes \K^{(\Lambda)})'
\end{equation}

  and similarly or by dualizing

\begin{equation}
\label{IsKL2}
  (\K^{\Lambda})' \hat{\otimes} (\K^{\Lambda})' \cong 
  (\K^{\Lambda} \hattimes \K^{\Lambda})'.
\end{equation}

\subsection*{Holomorphic vector-valued functions}
  The theory of holomorphic functions on an open subset $\Omega$ of $\C$ with
  values in a complete locally convex space $E$ over $\C$ parallels the 
  theory of holomorphic functions $g : \Omega \rightarrow \C$. The following 
  is easy to show.
\begin{Proposition}
  For a continuous $f : \Omega \rightarrow E$ the following properties
  are equivalent:
\begin{enumerate}
\item
  $f$ is holomorphic, i.e.~there is a continuous $f' : \Omega \rightarrow E$
  with 
\begin{equation}
  \lim \limits_{h \rightarrow 0} \, \frac{f (z + h) - f(z)}{h} \: = \: f' (z),
  \quad z \in \Omega.
\end{equation}
\item
  $f$ is analytic, i.e.~for each $z \in \Omega$ there are $c_n \in E$ such that
\begin{equation}
  f(z +h ) \: = \: \sum \limits_{n = 0}^{\infty} \, c_n \, h^n
\end{equation}
 for small $h \in \C$.
\item
  For all continuous linear forms $\alpha \in E'$ the function 
  $\alpha \circ f : \Omega \rightarrow \C$ is holomorphic.
\end{enumerate}
\end{Proposition}

Recall that for a given sequence $(c_n)$ of points in $E$ the power series 
$\sum \, c_n \, T^n$ converges in $E$ for $| h | < r$, if and only if
for all $ p \in {\rm cs} (E)$ and all $s < r$ the series 
$\sum \, p( c_n ) \, s^n$ converges (or equivalently:  is bounded).
Of course, the proposition has a straightforward generalization to open
subsets $\Omega$ of $\C^n$.

Let ${\cal O} ( \Omega , E)$ denote the space of holomorphic functions
$f : \: \Omega \rightarrow E$ endowed with the compact open topology.
Since $E$ is complete, ${\cal O} (\Omega , E)$ is  complete as well.
It can be shown that ${\cal O} (\Omega , E)$ induces on 
${\cal O} (\Omega) \otimes E$ the $\varepsilon$-topology
via the canonical map ${\cal O } \otimes E \ni f \otimes x 
\mapsto ( z \mapsto  f(z)\, x \in E) $. In the case of 
$\Omega \subset \C$ ( or $\Omega  \subset \C^n$) the nuclearity of 
${\cal O} (\Omega)$ implies that 
${\cal O} (\Omega) \hattimes E \cong {\cal O} (\Omega , E)$. 
Moreover, the space 
$ E_{\Omega} := {\cal O} (\Omega , E)$
is a ${\cal O} (\Omega)$-module with continuous action 
${\cal O} (\Omega) \times {\cal O} (\Omega , E ) \rightarrow
 {\cal O} (\Omega , E )$. 
\begin{Proposition}
  For complete locally convex spaces $E$ and $F$ the ${\cal O} (\Omega)$-module
  $(E \hattimes F)_{\Omega} = {\cal O} (\Omega , E \hattimes F )$ has the
  following universal property: Every continuous ${\cal O} (\Omega)$-bilinear
  map $\beta: E_{\Omega} \times F_{\Omega} \rightarrow G_{\Omega}$ can
  factored by a continuous ${\cal O} (\Omega)$-linear mapping 
  $\tilde{\beta} : \: (E \hattimes F)_{\Omega} \rightarrow G_{\Omega}$, 
  i.e.~$\beta = \tilde{\beta} \circ \otimes$.
\end{Proposition}
  Of course, any ${\cal O} (\Omega)$-linear $\lambda : \: 
  (E \hattimes F)_{\Omega} \rightarrow G_{\Omega}$ determines an 
  ${\cal O} (\Omega)$-bilinear mapping $\lambda \circ \otimes$ on 
  $E_{\Omega} \times F_{\Omega}$. Because of the property stated in the 
  proposition, we denote $(E \hattimes F)_{\Omega}$ also by $E_{\Omega}
  \hattimes_{{\cal O} (\Omega)} F_{\Omega}$. Indeed, $(E \hattimes F)_{\Omega}$
  is the completion of $E_{\Omega} \otimes_{{\cal O} (\Omega)} F_{\Omega}$.

\begin{Proof}
  $\beta$ is determined by $\beta |_{E \times F} : \: E \times F \rightarrow
  G_{\Omega} $. By the universal property of $E \hattimes F$ there is a 
  continuous $\C$-linear $\lambda : \: E \hattimes F \rightarrow G_{\Omega}$
  with $\beta |_{E \times F} = \lambda \circ \otimes$. $\lambda$ 
  determines a ${\cal O} (\Omega)$-bilinear map $\tilde{\beta}_0$ by
  $\tilde{\beta}_0 ( f \otimes (x \otimes y)) \: = \: f \cdot \lambda \,
  (x \otimes y)$
  for all $f \in {\cal O} (\Omega)$, $(x,y) \in E \times F$. Hereby ``$\cdot$''
  means pointwise multiplication of functions. $\tilde{\beta}_0$ is continuous
  on ${\cal O} (\Omega) \otimes (E \otimes F)$ and therefore has a continuous
  ${\cal O} (\Omega)$-linear extension $\tilde{\beta}$ to 
  ${\cal O} (\Omega) \hattimes ( E \hattimes F) = (E \hattimes F)_{\Omega}$
  with $\beta = \tilde{\beta} \circ \otimes$.
\end{Proof}
Another useful result for our presentation in section \ref{HoDeLoCoAl} is the 
following
\begin{Proposition}
  ${\cal L}_{{\cal O} (\Omega)} 
  (E_{\Omega} , F_{\Omega}) $ is algebraically and 
  topologically isomorphic to ${\cal O} (\Omega , {\cal L} (E,F))$. 
  Hereby ${\cal L}_{{\cal O} (\Omega)} (E_{\Omega} , F_{\Omega}) $ is the 
  ${\cal O} (\Omega)$-module of ${\cal O} (\Omega)$-linear and
  continuous maps $E_{\Omega} \rightarrow F_{\Omega}$ with the topology of 
  uniform convergence on the bounded sets of 
  $E_{\Omega} = {\cal O} (\Omega, E)$.
\end{Proposition}
\begin{Proof}
  We first show that  $T \mapsto T |_E$ defines an 
  isomorphism ${\cal L}_{{\cal O}
  (\Omega)} (E_{\Omega} , F_{\Omega} ) \rightarrow {\cal L} (E , F_{\Omega})$
  of locally convex spaces.
  Clearly, $t := T|_E : E \rightarrow F_{\Omega}$ is 
  $\C$-linear and continuous. Hence, $t \in {\cal L} (E, F_{\Omega})$,
  and the  map $T \mapsto t$ is ${\cal O} (\Omega)$-linear, 
  injective and continuous. For a given $s \in {\cal L} (E, F_{\Omega})$ 
\begin{equation}
  S (f \otimes x) \: := \: f \cdot s(x)
\end{equation}
  holds for all $f \in {\cal O} (\Omega)$, and $x \in E$ defines an element of 
  ${\cal L}_{{\cal O} (\Omega)} (E_{\Omega} , F_{\Omega} )$ with $S |_E = s$. 
  It remains to show that $t \mapsto T$ is continuous. 
  The topology on ${\cal L}_{{\cal O}(\Omega)} (E_{\Omega} , F_{\Omega} )$ 
  is generated by the seminorms
\begin{equation}
  q_{K,B}(T) := {\rm sup} \{q(Tf(z)) : f \in B, z \in K\}, \quad 
  T \in {\cal L}_{{\cal O}(\Omega)} (E_{\Omega} , F_{\Omega} ),
\end{equation}
  where $K$ is a compact subset of $\Omega$ and $B$ is a bounded subset of 
  $E_{\Omega}$. For such $K$ and $B$ the set 
  $B(K) := \{f(z) : f \in B \; {\rm and}  \; z \in K \} $ is a bounded   
  subset of $E$. Since $T$ is ${\cal O} (\Omega)$-linear one obtains 
  $T(\lambda \otimes x)(z) = \lambda(z)(Tx)(z) = t(\lambda(z)x)(z)$ whenever
  $z \in \Omega$, $\lambda \in {\cal O} (\Omega)$ and $x \in E$. Therefore, 
  $(Tf)(z) = t(f(z))(z)$ for all $z \in \Omega$. As a consequence,
\begin{eqnarray} 
\lefteqn{  q_{K,B}(T)  =  {\rm sup} \{q(t(f(z))(z)) : z \in K, f \in B \} 
  \leq } \nonumber \\
   & & \leq {\rm sup} \{q(tx(z)) : z \in K, f \in B(K) \} = q_{K,B(K)}(t),
\end{eqnarray}
  i.e. $t \mapsto T$ is continuous.
  
 	In order to prove ${\cal L}(E, F_{\Omega} ) \cong 
  {\cal O} (\Omega , {\cal L} (E,F))$, one can proceed in a similar manner.
  However, this isomorphism can also be deduced from general results on 
  the $\varepsilon$-product (cf. \cite{Sch:PrTe}) which is useful for 
  product formulas for spaces of holomorphic functions (cf. \cite{Sch:EPC}):
  Let $E_c'$ be the dual of the complete locally 
  convex space $E$ endowed with the topology of compact convergence. Then 
  $E \varepsilon F := {\cal L} (E_c',F)$, endowed with the topology of uniform 
  convergence on the equicontinuous subsets of $E'$, 
  is (by transposition) canonically isomorphic to 
  $F \varepsilon E$, and this new product is associative. By the theorem of  
  {\sc Mackey}, ${\cal L} (E,F) \cong E_c' \varepsilon F$, and for a  
  nuclear space $E$ one can show $E \varepsilon F \cong E \hattimes F$.
  Hence, ${\cal O} (\Omega,F) \cong {\cal O}(\Omega) \varepsilon F$ and we
  conclude
$$
{\cal L} (E , F_{\Omega} ) \: 
\cong \: E_c' {\varepsilon} \left( {\cal O} (\Omega ){\varepsilon} F \right) \: \cong \: {\cal O} (\Omega ) {\varepsilon} ( E_c' {\varepsilon} F ) \: 
\cong \: {\cal O} (\Omega) {\varepsilon} {\cal L} (E,F) \: 
\cong \: {\cal O} ( \Omega, {\cal L} (E,F) ).$$
  \end{Proof}

\begin{Corollary}
The following identifications hold for a complete locally convex space.
\begin{enumerate}
\item
  ${\cal L}_{{\cal O} (\Omega)} (E_{\Omega} , {\cal O} (\Omega)) \cong {\cal O} 
  (\Omega , E')$, or in a shorter way: $\left( E_{\Omega} \right)' \cong
  \left( E' \right)_{\Omega} = E'_{\Omega}$.
\item
  $E_{\Omega} \hattimes_{{\cal O} (\Omega)} E_{\Omega} 
  \cong (E \hattimes E)_{\Omega}$ and therefore
\begin{eqnarray}
  {\cal L}_{{\cal O} (\Omega)} (E_{\Omega} \hattimes_{{\cal O} (\Omega)} E_{\Omega} ,
  E_{\Omega} ) & \cong & {\cal O} (\Omega , {\cal L} (E \hattimes E, E) \\
  {\cal L}_{{\cal O} (\Omega)} ( E_{\Omega} , 
  E_{\Omega} \hattimes_{{\cal O} (\Omega)} E_{\Omega}) & \cong &
  {\cal O} (\Omega , {\cal L} (E, E \hattimes E).
\end{eqnarray}
\item
  If $E$ is strictly nuclear $E'_{\Omega} \hattimes_{{\cal O} (\Omega)}
  E'_{\Omega} \cong (E_{\Omega} \hattimes_{{\cal O} (\Omega)} E_{\Omega})'$.
\end{enumerate}
\end{Corollary}
\begin{Proof}
  The first two relations follow directly from the above proposition.
  The last isomorphism is a consequence of the following chain of
  isomorphisms:
\begin{equation}
  E'_{\Omega} \hattimes_{{\cal O} (\Omega)} E'_{\Omega} \: \cong \:
  {\cal O} (\Omega , E' \hattimes E' ) \: \cong \: {\cal O} (\Omega ,
  (E \hattimes E)')
\end{equation}
  since $E$ is strictly nuclear. Furthermore,
\begin{equation}
  {\cal O} (\Omega , (E \hattimes E)') \: \cong \: (E \hattimes E)'_{\Omega}
  \: \cong \: \left( (E \hattimes E )_{\Omega} \right)' \: \cong \: 
  \left( E_{\Omega} \hattimes_{{\cal O} (\Omega) } E_{\Omega} \right)'. 
\end{equation}
\end{Proof}

\nocite{BonFlaGerPin:HGSQGSDRPD, Mic:LMCTA, Sch:PrTe}
\bibliographystyle{amsplain}
\addcontentsline{toc}{section}{Bibliography}
\bibliography{pDGL,FuncAna,Algebra,ComAlgAna,MathPhys,TheoPhys,AlgGeo,Geo}

\ifx\undefined\bysame
\newcommand{\bysame}{\leavevmode\hbox to3em{\hrulefill}\,}
\fi
\begin{thebibliography}{10}

\bibitem{AleGroSch:CQHCSTI}
A.~Alekseev, H.~Grosse, and V.~Schomerus, {\em Combinatorial quantization of
  the {Hamiltonian} {Chern-Simons} theory {I}}, Comm.~Math.~Phys. {\bf 172}
  (1995), 317--358.

\bibitem{AleGroSch:CQHCSTII}
\bysame, {\em Combinatorial quantization of the {Hamiltonian} {Chern-Simons}
  theory {II}}, Comm.~Math.~Phys. {\bf 174} (1996), 561--604.

\bibitem{Ber:DLRT}
George~M. Bergman, {\em The diamond lemma for ring theory}, Advances in
  Mathematics {\bf 29} (1978), 178--218.

\bibitem{BonFlaGerPin:HGSQGSDRPD}
P.~Bonneau, M.~Flato, M.~Gerstenhaber, and G.~Pinczon, {\em The hidden group
  structure of quantum groups: Strong duality, rigidity and preferred
  deformations}, Comm.~Math.~Physics {\bf 161} (1994), 125--156.

\bibitem{ChaPre:GQG}
Vyjayanthi Chari and Andrew Pressley, {\em A {Guide} to {Quantum} {Groups}},
  Cambridge University Press, 1994.

\bibitem{Dri:QG}
V.~G. Drinfel'd, {\em Quantum groups}, Proc. of the ICM, 1986, pp.~798--820.

\bibitem{FadResTak:QLGLA}
L.~D. {Faddeev}, N.~Y. {Reshetikhin}, and L.~A. {Takhtajan}, {\em Quantization
  of {L}ie groups and {L}ie algebras}, Leningrad Math.~J. {\bf 1} (1990),
  193--225.

\bibitem{FocRos:PSMFCRSRM}
V.~V. Fock and A.~A. Rosly, {\em Poisson structure on moduli of flat
  connections on {Riemann} surfaces and {$r$}-matrix}, preprint ITEP {\bf
  72-92} (1992).

\bibitem{Ger:DRA}
Murray Gerstenhaber, {\em On the deformations of rings and algebras}, Ann. of
  Math. {\bf 79} (1964), 59--103.

\bibitem{Ger:DRAIV}
\bysame, {\em On the deformations of rings and algebras {$IV$}}, Ann. of Math.
  {\bf 99} (1974), 257--276.

\bibitem{GerSch:ACDT}
Murray Gerstenhaber and Samuel Schack, {\em Algebraic cohomology and
  deformation theory}, Deformation {Theory} of {Algebras} and {Structures} and
  {Applications}, Kluwer Academic publishers, 1988, pp.~11--264.

\bibitem{Gro:PTTEN}
Alexander Grothendieck, {\em Produits tensoriels topologiques et espace
  nucl\'eaires}, Memoirs of the AMS, vol.~16, Amer.~Math.~Soc., 1955.

\bibitem{Hor:ALPDO}
Lars {H\"ormander}, {\em The {Analysis} of {Linear} {Partial} {Differential}
  {Operators} {I-IV}}, Grund\-lehren, vol. 256(1983), 257(1983), 274(1985),
  275(1985), Springer-Verlag, 1983 - 1985.

\bibitem{Jim:QAUGLHAYBE}
M.~Jimbo, {\em A $q$-analogue of the {${U}_q (gl (N+1))$} {Hecke} algebra and
  the {Yang}-{Baxter} equation}, Lett.~Math.~Phys. {\bf 11} (1986), 247--252.

\bibitem{LarTow:QGQLA}
Richard~G. Larson and Jacob Towber, {\em Two dual classes of bialgebras related
  to the concepts of ``quantum group'' and ``quantum {Lie} algebra}, Comm.~in
  Alg. {\bf 19} (1991), 3295--3345.

\bibitem{Man:QGNG}
Yuri~I. Manin, {\em Quantum groups and non-commutative geometry}, Les
  Publications CRM, Universit\'e de Montr\'eal, 1988.

\bibitem{Mic:LMCTA}
Ernest~A. Michael, {\em Locally multiplicatively-convex topological algebras},
  Memoirs of the AMS, vol.~11, Amer.~Math.~Soc., 1952.

\bibitem{Pfl:LADQ}
Markus~J. Pflaum, {\em Local {Analysis} of {Deformation} {Quantization}}, Ph.D.
  thesis, {Fakult\"at} {f\"ur} Mathematik der
  {Ludwig-Maximilians-Universit\"at}, {M\"unchen}, November 1995.

\bibitem{Pfl:NCDQNOQCB}
\bysame, {\em A new concept of deformation quantization {I.} {Normal} order
  quantization on cotangent bundles}, Preprints des Sonderforschungsbereichs
  288, Berlin {\bf 186} (1995).

\bibitem{Pie:NLCS}
Albrecht Pietsch, {\em Nuclear {Locally} {Convex} {Spaces}}, Ergebnisse der
  Mathematik und ihrer Grenzgebiete, vol.~66, Springer-Verlag, Heidelberg,
  Berlin, New York, 1972.

\bibitem{Res:MuQu}
Nikolai Reshetikhin, {\em Multiparameter quantum groups and twisted
  quasi-triangular {H}opf algebras}, Lett. Math. Phys. {\bf 20} (1990),
  331--335.

\bibitem{ScheScho:MQMSFPB}
Peter Scheinost and Martin Schottenloher, {\em Metaplectic quantization of the
  moduli spaces of flat and parabolic bundles}, J.~reine angew.~Math. {\bf 466}
  (1995), 145--219.

\bibitem{Sch:EPC}
Martin Schottenloher, {\em $\epsilon$-product and continuation of holomorphic
  mappings}, Analyse fonctionelle et applications (Paris) (L.~Nachbin, ed.),
  Herman, 1975, pp.~261--270.

\bibitem{Sch:MiPrAlHoFu}
\bysame, {\em Michael problem and algebras of holomorphic functions},
  Arch.~Math. {\bf 37} (1981), 241--247.

\bibitem{Sch:PrTe}
Laurent Schwartz, {\em Produits tensoriels topologiques d'espaces vectoriels
  topologiques. {E}spaces vectoriels topologiques nucl\'eaires.
  {A}pplications}, S\'em.~1953/54, Inst.~Henri Poincar\'e.

\bibitem{Tak:HATAQGUSL}
M.~Takeuchi, {\em Hopf algebra techniques applied to the quantum group
  ${U}_q(sl (2))$}, Deformation {Theory} and {Quantum} {Groups} with
  {Applications} to {Mathematical} {Physics} (M.~Gerstenhaber and J.~Stasheff,
  eds.), Contemporary Mathematics, vol. 134, American Mathematical Society,
  1992, pp.~309--23.

\bibitem{Tak:STGL}
\bysame, {\em Some topics on ${{\rm {GL}}}_q (n)$}, J.~Algebra {\bf 147}
  (1992), 379--410.

\bibitem{Tak:LQG}
L.~A. Takhtajan, {\em Lectures on quantum groups}, Introduction to {Quantum}
  {Group} and {Integrable} {Massive} {Models} of {Quantum} {Field} {Theory}
  (Mo-Lin Ge and Bao-Heng Zhao, eds.), World Scientific, 1990.

\bibitem{Tre:TVSDK}
Tr\`eves, {\em Topological {Vector} {Spaces}, {Distributions} and {Kernels}},
  Academic Press Inc., New York, 1967.

\bibitem{WesZum:CDCQH}
J.~Wess and B.~Zumino, {\em Covariant differential calculus on the quantum
  hyperplane}, preprint CERN-TH-5697/90, 1990.

\end{thebibliography}
\end{document}